\documentclass[12pt]{article}
\usepackage{graphicx}
\usepackage{amsmath}
\usepackage{amssymb}
\usepackage{amsthm}
\usepackage{color}
\makeatletter
\renewcommand{\@makecaption}[2]{%
\vspace{\abovecaptionskip}%
\sbox{\@tempboxa}{#1. #2}
#1. #2\par
\vspace{\belowcaptionskip}}
\makeatother
\begin{document}
\begin{center}

\LARGE{Stability by linear approximation for time scale dynamical systems.}

\bigskip

{\large Sergey Kryzhevich${}^{a,b,}$\footnote{Corresponding author, email: kryzhevicz@gmail.com, s.kryzhevich@spbu.ru}, Alexander Nazarov${}^{a,c}$}

\bigskip

\textit{\noindent\small 
${}^a$ Faculty of Mathematics and Mechanics, Saint-Petersburg State University, Russia\\
${}^b$ School of Natural Sciences and Mathematics, The University of Texas at Dallas, TX, USA\\
${}^c$ Saint-Petersburg Department of Steklov Institute, Saint-Petersburg, Russia\\
}
\end{center}

\begin{abstract} {We study systems on time scales that are generalizations of classical
 differential or difference equations and appear in numerical methods. In this paper we consider linear systems 
 and their small nonlinear perturbations. In terms of time scales and of eigenvalues of matrices we
 formulate conditions, sufficient for stability by linear approximation.
For non-periodic time scales we use techniques of central upper Lyapunov exponents (a common tool 
of the theory of linear ODEs) to study stability of solutions. 
Also, time scale versions of the famous Chetaev's theorem on conditional instability are proved.
In a nutshell, we have developed a completely new technique in order to demonstrate that methods of non-authonomous linear ODE theory may work for time-scale dynamics.}

\end{abstract}


\section{Introduction}

We study dynamic equations on time scales i.e.~on unbounded closed subsets of $\mathbb R$.
 The time scale approach first introduced by S.\,Hilger and his collaborators (see [1]
 and references therein) was intensively developing during last decades. The first advantage of
 such approach is the common language that fits both for flows and diffeomorphisms. On the
 other hand, there are many numerical methods that correspond to non-uniform steps.
 Especially, this is applicable for modeling non-smooth or strongly non-linear dynamical systems.

Consider a motion of a particle in two distinct media, e.g. water and air. Evidently, to model
such system, it is not effective to use equidistant nodes. It is better to take more of them
inside time periods, corresponding to motions in water. This is a natural way to obtain a
non-trivial time scale in a real life problem (see [2] and references therein). Another
application of time scale analysis may be found for systems with delay [3]. 

In this sense it seems to be useful to generalize some results on stability theory, well-known 
for ODEs for time scale case. Mainly, we consider a general linear system (autonomous or
 non-autonomous) and its small non-linear perturbation. For the continuous dynamics, there
 exists a well-developed theory on stability by first approximation. For autonomous case
 there are classical stability criteria related to eigenvalues of a matrix of coefficients for
 linear approximation (call it ${\cal A}$).

For non-autonomous systems or for cases when eigenvalues do not give information on
 stability of the perturbed system,  there might be two approaches. The first one is based on
 the theory of Lyapunov functions (see [4] and references therein for review on the time scale
 version of this method). The second one involves integral inequalities, particularly the 
Gr\"{o}nwall--Bellmann inequality (see [5,6]). For ordinary differential
 equations there is a very powerful tool that allows to find stability of solutions via the so-called
 central upper Lyapunov exponents [7].  In this paper we combine all referred methods in order to study time scale dynamics.

Exponential stability for solutions of time-varying dynamic equations on a time scale have been
 investigated by many authors. We mention recent papers by Bohner and Martynyuk 
[8] (this article is also a good introduction to the theory of time scale systems), 
Du and Tien [9], Hoffacker and Tisdell [10] and Martynyuk [11]. 
 We also refer to papers [12--16] where related problems are studied and new approaches have been introduced.

A "multidimensional" analog of time scales called discrete differential geometry is also studied, see [17] and references therein. 
In such problems, time scales may appear, for instance, as discretizations of geodesic flows.

However, the following problems were open by now.

\begin{enumerate}
\item For constant matrices ${\cal A}$, are there any stability criteria for perturbed time-scale systems?
\item Is there any analog of Chetaev's theorem on instability by the first approximation for time scale systems?
\item Are there any sufficient conditions on stability by the first approximation, close to necessary ones?
\end{enumerate}

One of principal difficulties in the theory of time scale systems is that, generally speaking, in
 ``autonomous'' case (i.e. when the right hand side of the system does not depend on $t$) the
 system does not define a flow (a shift of a solution is not necesarily a solution, group
 property may be violated etc). Also, one must carefully check basic properties like
 smoothness of solutions that can be violated even for systems with smooth right hand sides.

In our paper we give positive answers to all mentioned questions. The main idea of our
 paper is very simple: methods of classical theory of linear non-autonomous differential
 equations are applicable for time scale systems. Here we notice that in time scale analysis 
 there are two types of derivatives: the so-called $\Delta$- and $\nabla$-derivatives (see [8] for details). 
 In this paper we study $\Delta$-derivatives only.
 However, it seems that the main ideas of our work can be easily transferred to equations with $\nabla$-derivatives.

We have two principal objectives. 

First, we provide sufficient conditions on stability by first approximation. We demonstrate
 that the obtained conditions are close to necessary ones. In our proofs, we use the techniques
 of central upper Lyapunov exponents. This approach seems to be novel for time scale
 analysis. We can use many related tools such as Millionschikov's rotations to obtain instability.
 
Secondly, we prove an analog of Chetaev's theorem on instability by first approximation.
 Specifics of time scales demands a novel, non-classical approach to proof since, generally
 speaking, we cannot use tools of the theory of autonomous systems, anymore. 

The paper is organized as follows. In Section 2 we give a brief introduction to time scale
 analysis. In Section 3, we give a review of existing results on stability of difference equations. 
 In these two sections we mostly refer to results of [8], changing some notations. In Section 4 we introduce one
 of the main ``characters'' of our paper: central upper exponent. Using this concept and the 
Gr\"{o}nwall--Belmann lemma, we study stability of solutions of time scale dynamical systems. 
In Section 5 we demonstrate that a time scale  system may be ``embedded'' into the ODE system. 
 In this connection we introduce so-called syndetic time scales i.e. time scales that do not have arbitrarily big gaps. We demonstrate that there is a correspondence 
 between linear time scale systems and linear systems of ordinary differential equations. We develop new technical tools that play a key role in proofs of following sections.
 In Sections 6 we provide a time scale generalization of the classical Millionschikov's result on
 attainability of the central upper exponent. The proof of Millionschikov's result is given in Appendix. As a corollary, we deduce a time scale version of
 the condition on instability by first approximation [18].   
The Lyapunov approach is developed in Section 7. Similarly to what happens in
 ODEs theory, we construct Lyapunov functions for time scale systems as quadratic forms
 and thus relate stability by first approximation with certain estimates on eigenvalues of the matrix of coefficients. 
 We also give a time scale version of Chetaev's instability theorem and a condition on instability by first approximation. 

\section{Time scale analysis}

We use following notation: $B(r,x)$ is the Euclidean ball in $\mathbb R^n$ with radius $r$, centered at $x$; $B_r=B(r,0)$, 
${\mathbf M}_{n}$ is the space of $n\times n$ complex
 matrices, $|\cdot |$ stands for a vector norm in ${\mathbb R}^n$ or ${\mathbb C}^n$
 and the corresponding operator norm. $E_n$ is $n\times n$ identity matrix.\medskip
 
\noindent\textbf{Definition 2.1.} A \emph{time scale} $\mathbb T$ is an unbounded closed subset of $\mathbb R$ with the inherited metric. 
Without loss of generality we assume that $0\in {\mathbb T}$.\medskip 

We consider two spaces of matrix functions: 
${\cal M}_R$ that is a space of continuous functions 
${\cal A}:{\mathbb R}\to {\mathbf M}_{n}$ and ${\cal M}_T$ that is the space of similarly defined functions ${\cal A}:{\mathbb T}\to {\mathbf M}_{n}$.

We set 
${\mathbb T}^+_a =[a,\infty) \bigcap {\mathbb T}$.\medskip

We introduce basic notions
 connected to the theory of time scales, which summarize the material from the recent book
 by Bohner and Peterson [19] (see also [20,21]).\medskip

\noindent\textbf{Definition 2.2.} Let $t \in {\mathbb T}$. We define 
the \emph{forward jump operator} $\sigma : {\mathbb T} \to {\mathbb T}$ by $\sigma (t) := \inf\{s \in {\mathbb T};s > t\}$.

If $\sigma(t) > t$, we say that $t$ is right-scattered, while if $\sigma(t) = t$, then $t$ is called right-dense.  
Denote by ${\mathbb S}$ the set of right-scattering points and by $\mathbb D$ the set of 
right-dense points. Evidently, ${\mathbb T}={\mathbb S}\bigcup {\mathbb D}$ is a disjoint union. We always assume that $\sup\mathbb S=+\infty$.

The \emph{graininess function} $\mu : {\mathbb T} \to [0, \infty)$ is defined by $\mu(t) := \sigma(t) - t$ (Fig.\, 1).

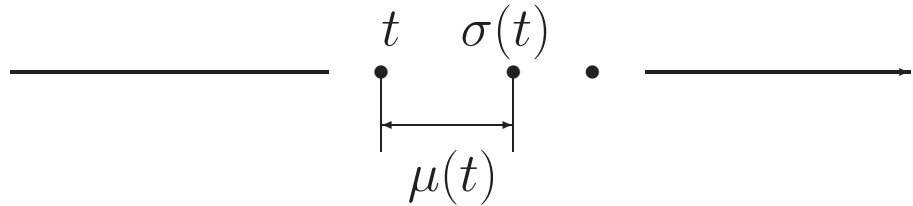
\begin{figure}[!ht]
\begin{center}
\begin{picture}(400,80)

\put(180,50){\circle*{5}}
\put(230,50){\circle*{5}}
\put(260,50){\circle*{5}}

\put(180,50){\line(0,-1){30}}
\put(230,50){\line(0,-1){30}}
\put(180,30){\vector(1,0){50}}
\put(230,30){\vector(-1,0){50}}

\linethickness{1pt}
\put(280,50){\vector(1,0){100}}
\put(40,50){\line(1,0){120}}

\put(190,5){\LARGE $\mu(t)$}
\put(210,60){\LARGE $\sigma(t)$}
\put(180,60){\LARGE $t$}

\end{picture}
\caption{A time scale.}
\end{center}
\end{figure}

\noindent\textbf{Definition 2.3.} A function $f : {\mathbb T} \to {\mathbb R}$ is called
 \emph{rd-continuous} provided it is continuous at right-dense points in ${\mathbb T}$
 and finite left-sided limits exist at left-dense points in ${\mathbb T}$. Denote the class of 
rd-continuous functions by ${\cal C}_{rd}  = {\cal C}_{rd}({\mathbb T}, {\mathbb R})$. We use the similar notation for vector and matrix functions.\medskip 

\noindent\textbf{Definition 2.4.} The function $f : {\mathbb T} \to {\mathbb R}$ is called 
$\Delta$-\emph{differentiable} at a point  $t \in {\mathbb T}$ if there exists 
$\gamma \in {\mathbb R}$ such that for any $\varepsilon > 0$ there exists a
 neighborhood $W$ of $t$ satisfying
$$|[f(\sigma(t)) - f(s)] - \gamma[\sigma(t) - s]| \le  \varepsilon|\sigma(t) - s|$$ 
for all $s \in W$. In this case we write $f^\Delta(t) = \gamma$. 

When 
${\mathbb T} = {\mathbb R}$, 
$x^\Delta(t) = \dot x(t)$. When  ${\mathbb T} = {\mathbb Z}$, $x^\Delta(n)$ is the standard forward difference operator 
$x(n + 1) - x(n)$.\medskip

\noindent\textbf{Definition 2.5.} If $F^\Delta(t) = f(t)$, $t \in {\mathbb T}$, then $F$ is a 
$\Delta$-\emph{antiderivative} of $f$, and the Cauchy $\Delta$-\emph{integral} is given by
$$\int_\tau^s f(t)\Delta t =F(s)-F(\tau), \qquad \mbox{for all} \quad s,\tau\in{\mathbb T}.$$

\noindent\textbf{Definition 2.6.} A function $p : {\mathbb T} \to {\mathbb R}$ is called \emph{regressive} provided that 
$1+\mu(t)p(t)\neq 0$ for all $t\in {\mathbb T}$ and \emph{positively regressive} if $1+\mu(t)p(t)> 0$ for all $t\in {\mathbb T}$. 
The set of all regressive and rd-continuous functions is denoted by 
${\cal R} = {\cal R}({\mathbb T}, {\mathbb R})$. The set of all positively regressive and rd-continuous function is denoted by ${\cal R}^+$.

A matrix mapping ${\cal A} : {\mathbb T} \to {\mathbf M}_n({\mathbb R})$ is called \emph{regressive} if for each 
$t \in {\mathbb T}$ the $n\times n$ matrix $E_n+\mu(t){\cal A}(t)$ is invertible, and \emph{uniformly regressive} if in addition the matrix function 
$(E_n+\mu(t){\cal A}(t))^{-1}$ is bounded.\medskip 

If ${\cal A}$ is constant, it is uniformly regressive if and only if
$$\inf\limits_{\mathbb T}|\lambda_k \mu(t)+1|> 0, \qquad k=1,\ldots,n \eqno (2.1)$$
where $\lambda_k$ are the eigenvalues of ${\cal A}$. Note that in this case solutions of the system
$$ x^\Delta={\cal A}x \eqno (2.2)$$
are unique and have finite Lyapunov exponents.

To prove this statement, it suffices to reduce $(2.2)$ to the normal form and thus reduce it to a set of linear first order equations.\medskip

\noindent\textbf{Definition 2.7.} For $p \in {\cal R}
$, we define \emph{the generalized exponential function} $e_p(t,s)$ 
by 
$$
e_p(t,s)=\exp \left(\int_s^t \xi_{\mu(\tau)}p(\tau)\, \Delta \tau \right),
\eqno (2.3)$$
where $\xi_h$ is the cylinder transformation given by formula
$$
\xi_h(z)=\begin{cases}\ \log(1+zh)/h  & \mbox{if} \quad h\neq 0;\\
\ z   & \mbox{if} \quad h = 0.
\end{cases}
$$ 
Note that the function $x(t) = e_p(t, t_0)x_0$ is the unique solution of the Cauchy problem
$$
x^\Delta(t) = p(t)x(t), \qquad x(t_0) = x_0, \qquad t_0 \in {\mathbb T},
$$
see [8] for details.\medskip

\noindent\textbf{Remark 2.8.} If $p\in{\cal R}$  then
$$1+ \int_a^t p(u)\Delta u\le e_p(t,a)\le \exp\left(\int_a^t p(u)\Delta u\right) \qquad \forall t\in{\mathbb T}^+_a.$$

\noindent\textbf{Theorem  2.9 (Comparison Theorem) [19].}
\emph{Let $t_0 \in {\mathbb T}$, $x, f \in {\cal C}_{rd}$ and $p \in {\cal R}^+$. Then 
$$x^\Delta(t) \le  p(t)x(t) + f(t), \qquad \mbox{for all} \quad t \in {\mathbb T}^+_{t_0}$$ 
implies 
$$x(t)\le x(t_0)e_p(t,t_0)+\int_{t_0}^t e_p(t,\sigma(\tau ))f(\tau )\Delta\tau,  \qquad 
 \forall \quad t\in {\mathbb T}^+_{t_0}.$$}

\noindent\textbf{Definition 2.10.} Let the matrix function ${\cal A}$ be regressive. The matrix function 
$\Phi_{\cal A}(t,t_0)$ satisfying 
$$\Phi^\Delta_{\cal A}(t,t_0) = {\cal A}(t)\Phi_{\cal A}(t,t_0), \qquad \Phi_{\cal A}(t_0, t_0) = E_n, \qquad t,t_0\in {\mathbb T}, \qquad t\ge t_0,$$
is called \emph{matrix exponential function} or \emph{fundamental matrix.} \medskip

\noindent\textbf{Theorem 2.11 [20].} \emph{Suppose that $n\times n$ matrix function ${\cal A}$ on the time scale 
$\mathbb T$ is regressive. Then  
\begin{itemize}
\item[$(i)$] the matrix exponential function $\Phi_{\cal A}(t,t_0)$ is uniquely defined for any $t_0\in {\mathbb T}$;  
\item[$(ii)$]  $\Phi_{\cal A}(t,r)\Phi_{\cal A}(r,s) =\Phi_{\cal A}(t,s)$ for $r,s,t \in {\mathbb T}$, $s\le r\le t$;
\item[$(iii)$] $\Phi_{\cal A}(\sigma(t), s) = (E_n + \mu(t){\cal A}(t))\Phi_{\cal A}(t, s);$
\item[$(iv)$] If ${\mathbb T} = {\mathbb R}$ and ${\cal A}$ is constant, then $\Phi_{\cal A}(t, s) = \exp({\cal A}(t-s))$;
\item[$(v)$] If ${\mathbb T}=h{\mathbb Z}$ with $h>0$ and ${\cal A}$ is constant, then $\Phi (t,s)=(E_n+h{\cal A})^\frac{t-s}h$.
\end{itemize}}

\section{Types of stability}

Let us consider a linear system
$$
x^\Delta={\cal A}(t) x 
\eqno (3.1) $$
and its nonlinear perturbation
$$
x^\Delta={\cal A}(t) x+ f(t,x). 
\eqno (3.2)$$
We always suppose that the matrix ${\cal A}(t)$ is bounded on $\mathbb T$.
When discuss systems $(3.1)$ and $(3.2)$, we denote by $x(t, t_0, x_0)$ the solution of the system subject to the initial value
$x(t_0)=x_0$.\medskip

\noindent\textbf{Definition 3.1.} 
\begin{itemize}
\item[a)] System $(3.1)$ is said to be \emph{stable} if, for every $t_0 \in {\mathbb T}$ and for every $\varepsilon > 0$ there exists 
a $\delta = \delta(\varepsilon,t_0) > 0$ such that
$$
|x_0| < \delta \qquad \Longrightarrow\qquad |x(t, t_0, x_0)| < \varepsilon,\quad \forall t\in T^+_{t_0}.
\eqno (3.3)$$
\item[b)] System $(3.1)$ is said to be \emph{uniformly stable} if for every $\varepsilon > 0$ there exists 
a $\delta = \delta(\varepsilon) > 0$ independent on initial point $t_0$, such that $(3.3)$ is satisfied.
\item[c)] System $(3.1)$ is said to be \emph{asymptotically stable} if it is stable and for any $t_0\in {\mathbb T}$ 
 there exists a positive value $c$ such that $|x_0| < c$ implies $x(t, t_0, x_0) \to 0$ as $t\to+\infty$.
\end{itemize}

\noindent\textbf{Theorem 3.2 (Choi and al. [22]; DaCunha [23]).} \emph{Let the matrix function ${\cal A}$ be regressive.
Linear system $(3.1)$ 
is stable if and only if all its solutions are bounded on 
${\mathbb T}^+_0$. It is uniformly stable if and only if there exists a positive constant $\gamma$, such that 
$|\Phi_{\cal A}(t,t_0)| \le  \gamma$, for all $t_0\in {\mathbb T}$, $t \in T^+_{t_0}$.}\medskip

Later on, when discuss systems $(3.1)$ and $(3.2)$ we always assume that the matrix function ${\cal A}$ is regressive. When discuss the 
system $(2.2)$ we always assume that ${\cal A}$ is uniformly regressive i.e. inequality $(2.1)$ is true unless the opposite statement is specified.

\section{Stability. Gr\"{o}nwall--Bellmann approach}

The principal objective of this section is to establish a condition on stability of a solution of a time scale system by first approximation. 
We use a tool well-known in the theory of linear systems. Namely, we introduce the so-called central upper exponents.\medskip  
 
\noindent\textbf{Definition 4.1.} A function $f:{\mathbb T}_0^+ \times B_r \to {\mathbb R}^n$ belongs to the class 
${\cal F}$ if it satisfies conditions
\begin{enumerate}
\item $f(t,0)=0$ for any $ t \geq 0;$
\item $\dfrac{\partial f}{\partial x}(t,0)=0$ for any $ t \geq 0;$
\item the Jacobi matrix 
$\dfrac{\partial f}{\partial x}(t,x)$ 
is uniformly continuous at ${\mathbb T}_0^+ \times B_r$.
\end{enumerate}

Observe that for any $f\in {\cal F}$ and any  $\varepsilon>0$ there exists $r_1>0$ such that
$$
|f(t,x)|\le \varepsilon |x|, \qquad \forall  t\in{\mathbb T}_0^+,\quad x\in B_{r_1}.
\eqno (4.1)$$
We consider systems $(3.1)$ and $(3.2)$ where $f\in {\cal F}$ for an $r>0$.\medskip

\noindent\textbf{Definition 4.2.} An rd-continuous function $u(t) : {\mathbb T}_0^+ \to {\mathbb R}$ is called \emph{upper function} 
for system $(3.1)$ if there exists a $C > 0$ such that the fundamental matrix of $(3.1)$
satisfies\footnote{See $(2.3)$ for the definition of generalized exponent $e_u(t,s)$.}
$$
| \Phi_{\cal A} (t,s)|  \le  C e_u(t,s),\qquad  \forall t,s \in {\mathbb T}_0^+,\quad t \ge s.
$$

Let $U$ be the set of all upper functions of $(3.1)$. We call the value
$$
\chi ({\cal A}): = \inf\limits_{u\in U} \chi_u: = \inf\limits_{u\in U} \limsup_{t \to \infty} \dfrac{1}{t} \int_0^t u(\tau)\Delta \tau
\eqno (4.2)$$
\emph{central upper exponent} for $(3.1)$.\medskip

Observe that this exponent is not less than the greatest Lyapunov exponent for solutions of system $(3.1)$. On the other hand, it is 
not greater than
$$\limsup_{t \to \infty} \dfrac{1}{t} \int_0^t |{\cal A}(s)|\, \Delta s.
$$

\noindent\textbf{Remark 4.3.} Central upper exponent may be greater than the greatest Lyapunov exponent. For linear systems of ODEs, corresponding 
example was given by Perron [24]. He controlled velocities of growth for solutions specifying the matrix of coefficients. Below in 
Example 4.6 we obtain the same effect for a system with a {\it constant} matrix and a special time scale 
 controlling switches between discrete and continuous regimes.\medskip

\noindent\textbf{Theorem 4.4.} \emph{If $\chi ({\cal A}) < 0$, there exists $\varepsilon > 0$ such that for any $f$ satisfying $(4.1)$, the zero 
solution of system $(3.2)$ is asymptotically stable.}\medskip

This statement is very close to Lemma 3.1 of [8].
To prove it we use the time scale version of Gr\"{o}nwall--Bellmann lemma first given in [19].\medskip 

\noindent\textbf{Lemma 4.5 (Gr\"{o}nwall--Bellmann Inequality).} \emph{Let $t_0 \in  {\mathbb T}$, $x, g, p\in  {\cal C}_{rd}$, $p \ge  0$. Then
$$x(t)\le g(t) + \int_{t_0}^t x(s) p(s) \, \Delta s \qquad \forall \quad t\in {\mathbb T}^+_{t_0}
$$
implies
$$x(t)\le g(t)+\int_{t_0}^t e_p(t,\sigma(s))g(s)p(s)\, \Delta s \qquad \mbox{for all } \quad t\in {\mathbb T}^+_{t_0}.
$$}

\noindent \textbf{Proof of Theorem 4.4.} Set $g(t,x)=f(t,x)\eta(|x|)$ where $\eta\in C^{\infty}(\mathbb R^+\to {\mathbb R})$ is a cut-off 
function such that $\eta(s)=1$ for $s\in (0,r_1/2)$, $\eta (s)=0$ for $s\ge r_1$, and $\eta'(s)\le0$ (here $r_1=r_1(\varepsilon)$ is a constant
from $(4.1)$). Evidently, the zero solution of the system $x^{\Delta}={\cal A}(t)x+g(t,x)$ is asymptotically stable if and only if one of $(3.2)$ is.

Hence, we may assume without loss of generality that inequality $(4.1)$ is satisfied for all 
$t\in {\mathbb T}_0^+$ and $x\in {\mathbb R^n}$. Then for any solution of $(3.2)$ we have
$$|x(t)|^\Delta=\dfrac{x^{\Delta}(t)\cdot x(t)}{|x(t)|}\le |{\cal A}(t)| |x(t)|+\varepsilon .
$$
By the Comparison Theorem, any solution of $(3.2)$ can be spread to ${\mathbb T}_0^+$.

Next, we rewrite $(3.2)$ in standard way as
$$x(t) = \Phi_{\cal A}(t,t_0)x(t_0) + \int_{t_0}^t \Phi_{\cal A}(t,s) f(s,x(s))\, \Delta s.
$$
By $(4.1)$, for any solution $x$ and for all $t \in {\mathbb T}^+_{t_0}$ we have
$$
|x(t)|\le |\Phi_{\cal A}(t,t_0) ||x(t_0)| + \varepsilon\int_{t_0}^t |\Phi_{\cal A}(t,s)|  |x(s)|\, \Delta s. 
\eqno (4.3)$$

Fix an upper function $u$ such that 
$\chi_u\in (\chi ({\cal A}), 0)$ where $\chi_u$ is defined by $(4.2)$. 
It follows from $(4.3)$ that
$$| x(t)| \le C | x(t_0 )| e_u(t,t_0)  + C \varepsilon \int_{t_0}^t e_u (t,s) |x(s)|\, \Delta s.
$$
Let $v(t)=|x(t)|/e_u (t,t_0)$. Then previous inequality can be rewritten as
$$v(t)\le C|x(t_0)|+C\varepsilon \int_{t_0}^t v(s) \Delta s.
$$
By the Gr\"{o}nwall--Bellmann inequality, we have
$$v(t)\le C |x(t_0)| + C^2\varepsilon |x(t_0)|\int_{t_0}^t e_{C\varepsilon}(t,\sigma(s))\, \Delta s.
$$
Taking into account Remark 2.8, we obtain
$$|x(t)|\le C |x(t_0)|\left(1+ C\varepsilon \int_{t_0}^t \exp(C\varepsilon(t-\sigma(s)))\, \Delta s \right) e_u (t,t_0).
$$
Evidently, if we choose $\varepsilon$ is sufficiently small i.e. $C\varepsilon<-\chi/2$, then  $x(t)$ tends to zero exponentially. 
This completes the proof. \hfill$\square$\medskip

Now we turn to the case of the constant matrix ${\cal A}$ and, respectively, to system $(2.2)$. Let $\lambda_k$, $k=1,\ldots, n$, be 
eigenvalues of the matrix ${\cal A}$. It is easy to see that the Lyapunov exponents are equal 
to\footnote{See Definition 2.7 for formula of the transformation $\xi_h$.}
$$
\nu_k =\limsup_{t\to \infty}  \dfrac 1t \int_{t_0}^t {\Re}\,\xi_{\mu(t)} (\lambda_k) \Delta t=\limsup_{t\to \infty}  \dfrac 1t\log|e_{\lambda_k}(t,t_0)|.
$$
To obtain this formula we proceed to the Jordan normal form in the corresponding system (2.2).
Then we apply Theorem 2.9 to estimate norms of solutions of obtained equations.

Evidently for any constant matrix $\cal A$ and any time scale $\mathbb T$ we have $\chi({\cal A})\ge\max\limits_k \nu_k$.
For ${\mathbb T}={\mathbb R}$ or ${\mathbb T}={\mathbb Z}$ the converse inequality is true, so $\chi({\cal A})=\max\limits_k \nu_k$.
However, in general this is not the case.\medskip

\noindent\textbf{Example 4.6.}
Take values $a\in (0,1)$ such that 
$$
\begin{array}{c}
\nu_1(a):=-2a + \log 11 (1-a)/6 <0;\\
\nu_2(a):= -0.5 a + \log 2 (1-a)/6<0;\\
\chi(a):= -0.5 a+\log 11 (1-a)/6>0.
\end{array}
\eqno (4.4) $$
(observe that the set of such values $a$ is open and non-empty, we can take $a$ such that $\chi(a)=0$ and slightly decrease it). We can choose $a$ rational, i.e. 
$a=p/q$, $p,q\in {\mathbb N}$.

Let ${\cal A}=\mathrm{diag}\, (-2,-0.5)$. 
Notice that for ${\mathbb T}_c={\mathbb R}$ the system $x^\Delta={\cal A}x$ is ''stable''. For the case 
${\mathbb T}_d=6{\mathbb Z}$ the situation is opposite (the Lyapunov exponents are $(\log 11)/6$ and $(\log 2)/6$ respectively).

For all $m,n\in {\mathbb N}$ consider the time scale 
$${\mathbb T}_{m,n}=\bigcup_{k\in {\mathbb Z}} [6qmk+6qn,6qmk+6pm+6qn] \bigcup 6 {\mathbb Z}.$$ 
This time scale is continuous for $a$-th part of the time and discrete for $(1-a)$-th part. By $(4.4)$, this implies $\nu_{1,2}<0$.

Now we take an increasing sequence $n_k\to \infty$ such that $m_k:=n_{k+1}-n_k\to \infty$ and $n_{k+1}/n_k\to 1$.

Introduce the time scale $\mathbb T$ by formula:
$${\mathbb T}\bigcap [6n_{j}q,6n_{j+1}q]={\mathbb T}_{m_j,n_j}\bigcap [6n_{j}q,6n_{j+1}q], \qquad j\in {\mathbb N}.$$

For this time scale, Lyapunov exponents of the system $x^\Delta={\cal A}x$ are still equal to $\nu_{1,2}$. Meanwhile, 
$\chi({\cal A})=\chi(a)>0$ that follows from the definition of central upper exponent, see also $(A.1)$.

\section{Reduction to ordinary differential equations}

The main objective of this section is to prove that a linear time scale system with a bounded and uniformly regressive matrix of coefficients 
can be ``embedded'' into a linear system of ordinary differential equations with a bounded matrix. Results of Lemmas 5.1 and 5.2 and Theorem 5.3 
given below are applied in the next section. However, they are of independent interest.\medskip

\noindent\textbf{Lemma 5.1.} \emph{For every linear regressive time scale system $(3.1)$, there exists a linear system of ordinary differential equations
$$
\dot x={\widetilde {\cal A}}(t)x,\qquad t\ge0,
\eqno (5.1)$$
such that if $\Phi_{\cal A}(t,s)$ is the fundamental matrix of system $(3.1)$ and ${\widetilde \Phi}_{\cal A}(t,s)$ is one for system $(5.1)$, then 
${\widetilde \Phi}_{\cal A}|_{{\mathbb T}^+_0\times {\mathbb T}^+_0}=\Phi_{\cal A}$.
If the matrix function ${\cal A}(t)$ is bounded and uniformly regressive, the matrix function ${\widetilde {\cal A}}(t)$ is bounded.}\medskip

\noindent\textbf{Proof.} 
We introduce the notation 
$$[t]_{\mathbb T}:=\sup ({\mathbb T}\cap (-\infty,t]).$$ 

We define $\widetilde{{\cal A}}(t)$ as follows:
$$
{\widetilde {\cal A}}(t)=\begin{cases}
{\cal A}(t), \qquad&\mbox{if} \quad t\in {\mathbb D}\cap \mathbb T^+_0;\\
\dfrac{1}{\mu([t]_{\mathbb T})} \mathrm{Log}\, (E_n+\mu([t]_{\mathbb T}){\cal A}([t]_{\mathbb T})), \qquad&
\mbox{if} \quad t\in {\mathbb R}_+\setminus{\mathbb D}.
\end{cases}
\eqno (5.2)$$
Notice that ${\widetilde {\cal A}}(t)$ is not uniquely defined. In general, it cannot be selected real even for a real matrix  ${\cal A}(t)$.
However, we can take it rd-continuous on $\mathbb R$ and constant on every connected subset of $\mathbb R\setminus{\mathbb T}$. Moreover, for any bounded and 
uniformly regressive matrix function ${\cal A}(t)$ we can select branches of matrix logarithm so that ${\widetilde {\cal A}}(t)$ is bounded. In what follows we fix 
such branches and use the notion $\log$ rather than $\mathrm {Log}$ that is the multivalued matrix logarithm.\medskip

Further, define ${\widehat\Phi}(t)$ by the formula
$$\widehat{\Phi}(t)=\begin{cases}
\Phi_{\cal A}(t,0), \qquad &\mbox{if}\quad t\in \mathbb T^+_0;\\
\exp({\widetilde {\cal A}}(t)(t-[t]_{\mathbb T}))\Phi_{\cal A}([t]_{\mathbb T},0), \qquad &\mbox{if}\quad
 t\in {\mathbb R}_+\setminus{\mathbb T}.
\end{cases}
$$ 
Direct calculation shows that $\widetilde{\Phi}_{\cal A}(t,s)={\widehat\Phi}(t){\widehat\Phi}^{-1}(s)$ 
is the fundamental matrix for system $(5.1)$. 
\hfill$\square$\medskip

The following statement is evident.\medskip

\noindent\textbf{Lemma 5.2.} \emph{Let ${\cal A}(t)$ be a bounded, rd-continuous and uniformly regressive matrix function on $\mathbb T^+_0$, and 
let ${\widetilde {\cal A}}(t)$ be defined by $(5.2)$. Then system $(3.1)$ is stable $($asymptotically stable$)$ 
if and only if corresponding system $(5.1)$ is.}\medskip

Given a bounded rd-continuous uniformly regressive matrix function ${\cal A}(t)$, we notice that for arbitrary rd-continuous matrix function ${\cal B}(t)$ 
such that $\|{\cal B}\|_{{\mathbb L}^\infty}\le \delta$, $({\cal A}+{\cal B})(t)$ is also uniformly regressive provided that $\delta>0$ is sufficiently small.
Thus, we can define the extension $(\widetilde {{\cal A}+{\cal B}})(t)$ similarly to $(5.2)$, and it is easy to see that the nonlinear mapping
$$
{\cal L}[{\cal B}]\equiv\widehat {\cal B}:=\widetilde {{\cal A}+{\cal B}}-{\widetilde {\cal A}}
$$
is Lipschitz.\medskip

Now let $\widetilde\Phi(t,s)$ be the fundamental matrix for the system
$$
\dot x=(\widetilde {\cal A}(t)+\widehat {\cal B}(t))x.
\eqno (5.3)$$
We define
$$
{\cal B}(t)=:\widetilde{\cal L}[\widehat {\cal B}]=\begin{cases}
\widehat {\cal B}(t), \qquad & \mbox{if} \qquad t\in {\mathbb D}\cap \mathbb T^+_0;\\
\dfrac{\widetilde\Phi(\sigma(t),t)-E_n}{\mu(t)}-{\cal A}(t), \qquad &
\mbox{if} \qquad t\in {\mathbb S}\cap \mathbb T^+_0.
\end{cases}
\eqno (5.4) $$

We are in the position to prove the main result of this section.\medskip

\noindent\textbf{Theorem 5.3.} \emph{Suppose that ${\cal A}(t)$ is a bounded, rd-continuous and uniformly regressive matrix function on $\mathbb T^+_0$, and 
let ${\widetilde {\cal A}}(t)$ be defined by $(5.2)$. Let a continuous matrix function $\widehat {\cal B}(t)$ on $\mathbb R_+$ satisfy the following assumptions:
$$\|\widehat {\cal B}\|_{{\mathbb L}^\infty}\le \delta;$$
$$\widehat {\cal B} (t)=0  \qquad \mbox{for all} \quad t\in {\mathbb R},\quad  \mathrm{dist}\, (t,{\mathbb T})>S,
\eqno (5.5) $$
where $S>0$ is arbitrary given constant while $\delta>0$ is sufficiently small. Then
formula $(5.4)$ defines the rd-continuous matrix function ${\cal B}(t)$ on $\mathbb T^+_0$ such that $({\cal A}+{\cal B})(t)$ is uniformly regressive.
Moreover, the nonlinear operator $\widetilde{\cal L}$ defined by $(5.4)$ is Lipschitz left inverse to  $\cal L$, and its Lipschitz constant depends
only on ${\cal A}$, $S$ and $\delta$ $($in particular, it does not depend on the time scale $\mathbb T${}$)$.}\medskip

\noindent\textbf{Proof.} Fix a matrix ${\widehat {\cal B}}_0$ subject to $(5.5)$, such that $\|\widehat {\cal B}\|_{{\mathbb L}^\infty}=1$, and define 
$\Phi_\delta (t,s)$, $\delta\in[0,1]$, as the Cauchy
matrix for the system
$$\dot x=({\widetilde {\cal A}}(t)+\delta {\widehat {\cal B}}_0(t))x.
$$

Let $R_\delta:=\widetilde{\cal L}(\delta {\widehat {\cal B}}_0)$. By $(5.4)$, 
$$\dfrac{\partial}{\partial \delta}R_\delta(t)=\begin{cases}
\widehat {\cal B}_0(t)\vphantom{\dfrac{1}{\mu(t)}}, \qquad & \mbox{if} \qquad t\in {\mathbb D}\cap \mathbb T^+_0;\\
\dfrac{1}{\mu(t)} \dfrac{\partial \Phi_\delta(\sigma(t),t)}{\partial \delta}, \qquad &
\mbox{if} \qquad t\in {\mathbb S}\cap \mathbb T^+_0.
                  \end{cases}
$$

It is easy to see that 
$\partial \Phi_\delta(\sigma(t),t)/\partial \delta=U(\sigma(t))$, where $U(s)$ is the solution of
the matrix initial value problem
$$
\dot U(s)=({\widetilde {\cal A}}(s)+\delta {\widehat {\cal B}}_0 (s))U(s)+{\widehat {\cal B}}_0 (s) \Phi_\delta(s,t); \qquad U(t)=0.
$$
Consequently,
$$
\dfrac{\partial \Phi_\delta(\sigma(t),t)}{\partial \delta}=\int_{t}^{\sigma(t)} \Phi_\delta (\sigma(t),s) 
{\widehat {\cal B}}_0(s)\Phi_\delta(s,t)\, ds.
\eqno (5.6) $$

Without loss of generality we can assume that $S\ge1$. If $0<\mu(t)\le 2S$ then
\begin{eqnarray*}
&\left|\dfrac{\partial \Phi_\delta(\sigma(t),t)}{\partial \delta}\right|
&\le\int_t^{\sigma(t)} |\Phi_\delta (\sigma(t),s) |\, |{\widehat {\cal B}}_0(s)|\,|\Phi_\delta(s,t)|\, ds\\
&&\le \int_t^{\sigma(t)}\exp((\|\widetilde {\cal A}\|_{{\mathbb L}^\infty}+\delta)(\sigma(t)-s))\exp((\|\widetilde {\cal A}\|_{{\mathbb L}^\infty}+\delta)(s-t))\,ds\\
&&\le \mu(t)\exp((\|\widetilde {\cal A}\|_{{\mathbb L}^\infty}+\delta)\mu(t)),
\end{eqnarray*}
and thus
$$
\left|\dfrac{\partial}{\partial \delta}R_\delta(t)\right|\le\exp(2S(\|\widetilde {\cal A}\|_{{\mathbb L}^\infty}+1))=:M.
$$
Otherwise, for $\mu(t)>2S$, from $(5.6)$ we obtain the following estimate for the function
$w(\delta)=|\Phi_\delta (\sigma(t),t) |$:
\begin{eqnarray*}
&\dfrac{dw(\delta)}{d\delta}&\le\left|\dfrac{\partial \Phi_\delta(\sigma(t),t)}{\partial \delta}\right|\le
\int_t^{\sigma(t)} |\Phi_\delta (\sigma(t),s) |\, |{\widehat {\cal B}}_0(s)|\,|\Phi_\delta(s,t)|\, ds\\
&&\le\int_{t}^{t+S}
|\Phi_\delta (\sigma(t),t) |\,|\Phi_\delta (t,s) | \,|{\widehat {\cal B}}_0(s)|\,|\Phi_\delta(s,t)|\, ds\\
&&+\int_{\sigma(t)-S}^{\sigma(t)}|\Phi_\delta (\sigma(t),s) |\,|{\widehat {\cal B}}_0(s)|\,|\Phi_\delta(s,\sigma(t))|\,|\Phi_\delta (\sigma(t),t) |\, ds\\
&&\le 2S \big(\exp((\|\widetilde {\cal A}\|_{{\mathbb L}^\infty}+\delta)S)\big)^2 w(\delta)\le 2SM w(\delta).
\end{eqnarray*}
Since $w(0)=|E_n+\mu(t){\cal A}(t)|$, we obtain
$$\left|\dfrac{\partial}{\partial \delta}R_\delta(t)\right|\le 2SM\exp(2\delta SM)\cdot\frac{|E_n+\mu(t){\cal A}(t)|}{\mu(t)},
$$
and the derivative of $R_\delta$ is uniformly bounded for all $t$.\medskip

The identity $\widetilde{\cal L}[{\cal L}[{\cal B}]]\equiv {\cal B}$ is established by direct calculation. Finally, for sufficiently small $\delta$
the matrix function $({\cal A}+{\cal B})(t)$ is evidently uniformly regressive.\hfill$\square$\medskip

\noindent\textbf{Definition 5.4.} We say that the time scale $\mathbb T$ is \emph{syndetic}, if $\limsup\limits_{{\mathbb T}\ni t\to+\infty} \mu(t)<+\infty$.
This notion is similar to one used in Combinatorics and Number Theory.\medskip

\noindent\textbf{Corollary 5.5.} \emph{Let $\mathbb T$ be syndetic time scale. Then we can select the value $S$ so that the relation $(5.5)$ is satisfied for all 
matrix functions $\widehat {\cal B}$.}\medskip

\noindent\textbf{Remark 5.6.} For non-syndetic time scales operator $\widetilde{\cal L}$ is not continuous without assumption $(5.5)$. Namely, there exists a small 
matrix $\widehat {\cal B}$ such that $\widetilde{\cal L} [\widehat {\cal B}]$ is unbounded with respect to $t$.

\section{Instability via Millionschikov's rotations}

Now we prove two statements in a sense converse to Theorem 4.4.\medskip

{
\noindent\textbf{Theorem 6.1.} \emph{Let the matrix ${\cal A}(t)$ in $(3.1)$ be bounded and uniformly regressive, and let $\chi({\cal A})\ge 0$. Suppose that the time 
scale $\mathbb T$ is syndetic. Then for any $\delta>0$ there exists a matrix ${\cal B}_\delta (t)$ such that 
$$
\|{\cal B}_\delta \|_{{\mathbb L}^\infty}\le \delta,
\eqno (6.1) $$
and the system
$$
x^\Delta=({\cal A}(t)+{\cal B}_\delta (t)) x
\eqno (6.2)$$
is unstable.} \medskip

We reproduce the classical result on attainability of upper center exponents for systems of ordinary differential equations. 
In Appendix we provide a full proof of that result cause we need its details later on. We slightly modify the original proof in order to be able to apply Theorem 5.3.\medskip

\noindent\textbf{Theorem 6.2 (Millionschikov) [25].} 
\emph{Consider a system $(5.1)$ and suppose that $\|\widetilde {\cal A}\|_{{\mathbb L}^\infty}=a<\infty$. 
Then for any $\varepsilon>0$ and $\delta>0$,  there exists a continuous matrix $\widehat {\cal B}(t)$ such that 
$$
\|\widehat {\cal B}\|_{{\mathbb L}^\infty}\le (2a+1) \delta \eqno (6.3)
$$
and the greatest Lyapunov exponent of system $(5.3)$
is greater than $\chi(\widetilde {\cal A})-\varepsilon$, where $\chi(\widetilde {\cal A})$ is the central upper exponent of system $(5.1)$.} \medskip

\noindent\textbf{Proof of Theorem 6.1.} 

\noindent\textbf{1.} First, we assume $\chi({\cal A})>0$.

We embed the system $(3.1)$ to a linear system of ODEs $(5.1)$ and observe the evident fact that $\chi(\widetilde {\cal A})\ge\chi({\cal A})$. Denote 
$a:=\|\widetilde {\cal A}\|_{{\mathbb L}^\infty}$. 

By Theorem 6.2 and Corollary A.1, for any $\delta_1>0$ there exists a continuous perturbation $\widehat {\cal B}(t)$ such that 
$\|\widehat {\cal B}\|_{{\mathbb L}^\infty}\le (2a+1) \delta_1$, and the greatest Lyapunov exponent of ODE system $(5.3)$ is greater than $\chi({\cal A})/2$.
Since the time scale is syndetic, the reduction of $(5.3)$ to ${\mathbb T}$ does also have a positive Lyapunov exponent. 

Denote ${\cal B}_\delta=\widetilde{\cal L} [\widehat {\cal B}]$. By Theorem 5.3 (with regard to Corollary 5.5), there exists a $K>0$ such that
$\|{\cal B}_\delta\|_{{\mathbb L}^\infty}<K(2a+1) \delta_1$ for small values of $\delta_1$.  

To finish the proof it suffices to take $\delta_1\le \frac {\delta}{K(2a+1)}$.\medskip

\noindent\textbf{2.} Now we study the case $\chi ({\cal A})=0$. First of all, observe that the transformation 
$y=\exp (\varepsilon t)x$ transfers a system $\dot x=P(t)x$ to $\dot y=(P(t)+\varepsilon E_n)y$.

We begin with the same procedure as in part {\bf 1} and construct a perturbation $\widehat {\cal B}'(t)$ such that 
$\|\widehat {\cal B}'\|_{{\mathbb L}^\infty}\le \frac{\delta}{3K(2a+1)}$, and the greatest Lyapunov exponent of the system
$$\dot x=(\widetilde {\cal A}(t)+\widehat {\cal B}'(t))x$$
is greater than $-\frac {\delta}{3K(2a+1)}$.

Now it is easy to see that the matrix ${\cal B}_\delta=\widetilde{\cal L} \big[\widehat {\cal B}'+\frac {2\delta}{3K(2a+1)}E_n\big]$ satisfies $(6.1)$
and provides unstable system  $(6.2)$.\hfill$\square$\medskip

For non-syndetic time scales the similar result is true under the positivity assumption for the upper central exponent.\medskip

\noindent\textbf{Theorem 6.3.} \emph{Let the matrix ${\cal A}(t)$ in $(3.1)$ be bounded and uniformly regressive, and let $\chi({\cal A})> 0$.
Then for any $\delta>0$ there exists a matrix ${\cal B}_\delta (t)$ satisfying $(6.1)$ and such that system $(6.2)$ is unstable.} \medskip

\noindent\textbf{Proof.} Notice that direct repetition of the proof of Theorem 6.1 does not work since the assumption $(5.5)$ for the perturbation
$\widehat {\cal B}(t)$ constructed in Theorem 6.2 may be violated. So, we need to modify the Millionschikov method.\medskip

As in the proof of Theorem 6.1, we embed system $(3.1)$ to linear system of ODEs $(5.1)$ and denote $a:=\|\widetilde {\cal A}\|_{{\mathbb L}^\infty}$. 
Recall that 
$${\widetilde {\cal A}}(t)=\dfrac{1}{\mu([t]_{\mathbb T})} \log\, (E_n+\mu([t]_{\mathbb T}){\cal A}([t]_{\mathbb T})),\qquad t\in {\mathbb R}_+\setminus{\mathbb D}.
$$
We choose $S_1(\varepsilon, a)$ such that $\log\, (1+sa)/s<\varepsilon/4$ provided $s>S_1$.\medskip

Now we follow the proof of Theorem 6.2 (see Appendix) up to definition of segments $Q_j$ (Fig.~4). Without loss of generality we assume 
that $T_0\ge S_1\max(1,8a/\varepsilon)$.

We start with the segment $[T,2T]$ where $T$ is defined in Step 2 of the proof of Millionschikov's theorem, see $(A.6)$. 

There may be a segment $[\tau_{j},\tau_{j}+1]$ and a segment $[\tau_{j+1},\tau_{j+1}+1]$ where the perturbation from Theorem 6.2 is non-zero. On these segments 
two steps of rotation from $x_0(t)$ to $x_1(t)$ are performed.

Observe that
$$[\tau_j,\tau_{j+1}]\bigcap {\mathbb T}\neq \emptyset. \eqno (6.4)
$$
Otherwise, the segment $[\tau_j,\tau_{j+1}]$ completely belongs to a "gap"\ of the time scale of length greater than $S_1$ and, consequently, inequality $(A.8)$ fails 
for $i=j$.  

We introduce time instants $\tau_j'$ and $\tau_{j+1}'$ as follows.
$$
\tau_j'=\begin{cases}
         \tau_j, & \mbox{if}\quad \mathrm{dist}\, (\tau_j, {\mathbb T})\le S_1;\\
         \sigma([\tau_j]_{\mathbb T}), & \mbox{if}\quad \mathrm{dist}\, (\tau_j, {\mathbb T})> S_1;
        \end{cases}
$$
$$
\tau_{j+1}'=\begin{cases}
         \tau_{j+1}, & \mbox{if}\quad \mathrm{dist}\, (\tau_{j+1}, {\mathbb T})\le S_1;\\
         [\tau_{j+1}]_{\mathbb T}, & \mbox{if}\quad \mathrm{dist}\, (\tau_{j+1}, {\mathbb T})> S_1.
        \end{cases}
$$
Notice that $\tau_{j+1}\ge\tau_{j+1}'\ge \tau_j'\ge\tau_j$ by virtue of $(6.4)$.

By definition of $S_1$, the following analog of inequality $(A.8)$ is satisfied in any case:
$$
\dfrac{| x_1(\tau'_{j+1})|}{| x_1(\tau'_j)|}:\dfrac{| x_0(\tau'_{j+1})|}{| x_0(\tau_j')|} \geq
\exp\left(\dfrac{\varepsilon T_0}{2}-\dfrac{\varepsilon T_0}4\right)=\exp\left(\dfrac{\varepsilon T_0}{4}\right).
\eqno (6.5)$$
This implies $\tau_j'+1<\tau_{j+1}'$ (recall that $aS_1\le \varepsilon T_0/8$). Consequently,
$$[\tau'_j, \tau'_j+1]\bigcap [\tau'_{j+1}, \tau'_{j+1}+1]=\emptyset.
$$

Observe that inequality $(6.5)$ is still enough to imply item \textbf{C} of the Step 3 in the proof of Theorem 6.2 (see also the footnote to Eq.
$(A.3)$).

So, the proof of Theorem 6.2 still passes with $\tau_j$ and $\tau_{j+1}$ replaced with $\tau_j'$ and $\tau_{j+1}'$. On the interval $[T,2T]$ the perturbation 
$\widehat {\cal B} (t)$ is non-zero on segments $[\tau'_j, \tau'_j+1]$ and $[\tau'_{j+1}, \tau'_{j+1}+1]$ only. The similar statement is true for all other 
segments $[kT,(k+1)T]$.

Thus, for any $\delta_1>0$ and $\varepsilon>0$ we have constructed a continuous perturbation $\widehat {\cal B}(t)$ such that 
$\|\widehat {\cal B}\|_{{\mathbb L}^\infty}\le (2a+1) \delta_1$, the greatest Lyapunov exponent of ODE system $(5.3)$ is greater than $\chi({\cal A})-\varepsilon$,
and the inequality $(5.5)$ holds with $S=S_1+1$.

Denote ${\cal B}_\delta=\widetilde{\cal L} [\widehat {\cal B}]$. By Theorem 5.3, there exists a $K(\varepsilon,a)>0$ such that
$\|{\cal B}_\delta\|_{{\mathbb L}^\infty}<K(2a+1) \delta_1$ for small values of $\delta_1$. 

Now we put $\varepsilon=\chi({\cal A})/2$ and claim that the greatest Lyapunov exponent of the time scale system $(6.2)$ is not less than $\chi({\cal A})/2$. Indeed, 
let $\Phi(t)$ be a fundamental matrix of system $(5.3)$. By construction, there exists an unbounded sequence $t_k\in {\mathbb R}$ such that 
$|\Phi(t_k)|>\exp (\chi({\cal A})t_k/2)$. 

Denote $s_k=[t_k]_{\mathbb T}$ and notice that if $t_k-s_k>S_1$ then $\mu(s_k)>S_1$ and therefore $|{\widetilde {\cal A}}(t)|\le\varepsilon/4$
for $t\in[s_k,t_k]$. This gives
\begin{eqnarray*}
|\Phi(s_k)|&\ge& |\Phi(t_k)|\min\{\exp(-\varepsilon(t_k-s_k)/4), \exp(-aS_1)\}\\
&\ge&\exp(\chi({\cal A}) t_k/2-\varepsilon (t_k-s_k)/4) \exp(- aS_1)\\
&\ge&\exp(-aS_1) \exp (\chi({\cal A})s_k/2),
\end{eqnarray*}
and the claim follows.

To finish the proof it suffices to take $\delta_1\le \frac {\delta}{K(2a+1)}$.\hfill$\square$\medskip

The next statement is a generalization of the result of [18] for time scales.\medskip

\noindent\textbf{Theorem 6.4.} \emph{Let the matrix ${\cal A}(t)$ in $(3.1)$ be bounded and uniformly regressive, and let $\chi({\cal A})> 0$. 
Then there exists a continuous map $f:{\mathbb T}_0^+ \times B_1 \to {\mathbb R}^n$ such that\footnote{See Definition 4.1.} $f\in{\cal F}$ and 
the solution $x(t)\equiv0$ of the corresponding system $(3.2)$ is unstable.
If the time scale $\mathbb T$ is syndetic, the same is true provided $\chi({\cal A})\ge 0$.}\medskip

\noindent\textbf{Proof.}  By Theorems 6.1 and 6.3, there exist continuous perturbation matrices ${\cal B}_{\ell}(t)$, $\ell\in {\mathbb N}$, such that 
\begin{itemize}
\item[1)] $|{\cal B}_{\ell}(t)| \le 2^{-{\ell}}$ for all ${\ell}\in {\mathbb N}$, $t\in {\mathbb T}$;
\item[2)] for every ${\ell}\in {\mathbb N}$ the system
$$
x^\Delta=({\cal A}(t)+{\cal B}_{\ell}(t))x
\eqno (6.6)$$  
has a solution with positive Lyapunov exponent.
\end{itemize}

Fix an unbounded solution $x_1(t)$ of system $(6.6)$ for ${\ell}=1$ such that $|x_1(0)|\le1$. Select $T_1>0$ so that 
$|x_1(T_1)|\ge 2,$ $|x_1(t)|<2$ while $0 \leq t<T_1.$ Then we construct an unbounded solution $x_2(t)$ of system  $(6.6)$ for ${\ell}=2$ such that 
$|x_2(t)|<|x_1(t)|/2$ for $0 \leq t \leq T_1.$ Given $x_2(t)$ we select the first time instant $T_2$ such that $|x_2(T_2)|\ge 2$. Then we 
construct $x_3(t)$ and $T_3$ and so on (Fig.\, 2).\medskip

\begin{figure}[!ht]
\begin{center}
\includegraphics*[width=4in]{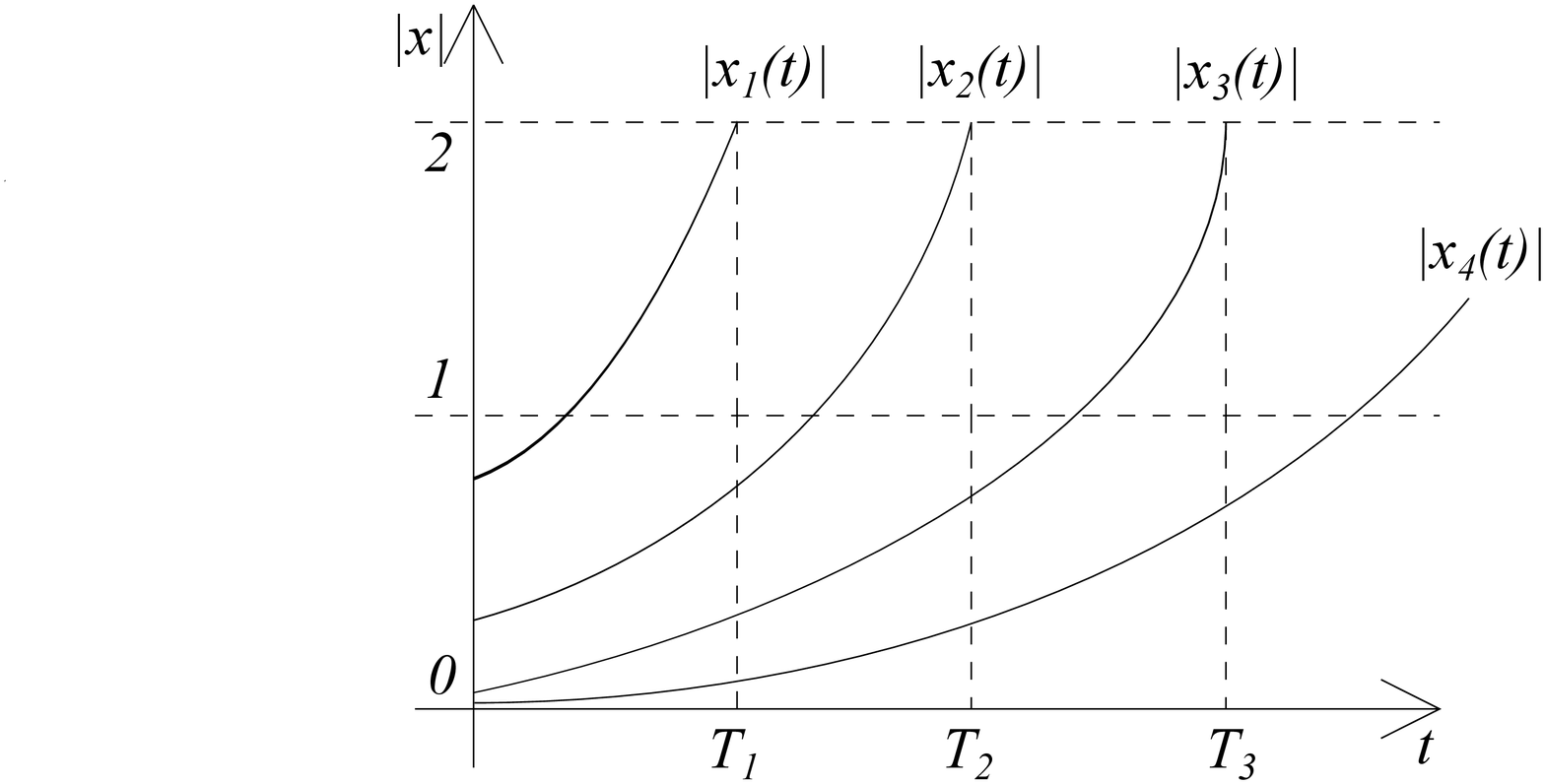}
\end{center}
\caption{Solutions of the perturbed system (${\mathbb T}={\mathbb R}$).}
\end{figure}

Now we construct a map $f: [0,\infty) \times B_1 \to {\mathbb R}^n$ such that the system $(3.2)$ coincides with $(6.6)$ in some neighborhood of the graph of
$x_{\ell}$ on $[0,T_{\ell}]$ for all ${\ell}\in{\mathbb N}$. This implies that all $x_{\ell}(t)$ are solutions of $(3.2)$ on $[0,T_{\ell}]$. 
Since $|x_{\ell}(0)| \to 0$ as ${\ell}\to\infty$ and $|x_{\ell}(T_{\ell})|\ge 2$, the zero solution of $(3.2)$ is unstable.\medskip

We set%
$$\Psi_{\ell}(t,x)=\begin{cases}
{\cal B}_{\ell}(t)\phi\big(\frac {|x-x_{\ell}(t)|}{\epsilon_{\ell}(t)}\big), &\mbox{if\quad} |x-x_{\ell}(t)| \leq \epsilon_{\ell}(t),\ \ t<T_{\ell}; \\
0,&\mbox{otherwise,}
                   \end{cases}
$$
where $\phi(s)$ is a smooth cut-off function equal to one for $s\le 1/2$ and to zero for $s\ge1$, while $\epsilon_{\ell}(t)=| x_{\ell}(t)|/4$. 
Evidently, all $\Psi_{\ell}(t,x)$ are continuous and continuously differentiable with respect to $x$.

We define
$$
f(t,x)=\sum_{{\ell}=1}^{\infty} \Psi_{\ell}(t,x)\cdot x.
$$
Notice that for any $k<{\ell}$ and $t<T_k$ we have $|x_{\ell}(t)|<\frac {|x_{{\ell}-1}(t)|}2<\dots <\frac {|x_k(t)|}{2^{{\ell}-k}}$
and therefore $|x_{\ell}(t)-x_k(t)|>|\epsilon_{\ell}(t)|+|\epsilon_k(t)|$. Thus, balls $B(x_k(t),\epsilon_k(t))$ and $B(x_{\ell}(t),\epsilon_{\ell}(t))$ 
are pairwise disjoint, and the map $f$ satisfies the above assumption. Moreover, evidently $f(t,0)\equiv0$. So, we should only examine
differentiated series
$$
\sum_{{\ell}=1}^{\infty} \Psi_{\ell}(t,x)+\sum_{{\ell}=1}^{\infty}\nabla_x\Psi_{\ell}(t,x)\cdot x.
\eqno (6.7)$$

Since $|\Psi_{\ell}(t,x)| \le|{\cal B}_{\ell}(t)| <1/2^{\ell}$, the first series in $(6.7)$ uniformly converges on ${\mathbb T}^+_0\times B_1$.
Further, for $(t,x)$ in support of $\Psi_{\ell}$ we have
\begin{eqnarray*}
&&|\nabla_x\Psi_{\ell}(t,x)\cdot x|\le |{\cal B}_{\ell}(t)|\,\dfrac{\max_s|\phi'(s)|}{\epsilon_{\ell}(t)}\,|x|\\
&&\le \dfrac{\max_s|\phi'(s)|}{2^{\ell}}\dfrac{| x_{\ell}(t)|+| x-x_{\ell}(t)|}{\epsilon_{\ell}(t)}\le\dfrac{5\max_s|\phi'(s)|}{2^{\ell}},
\end{eqnarray*}
and the second series in $(6.7)$ also uniformly converges. Therefore, the sum of the series $(6.7)$ equals $\dfrac{\partial f}{\partial x}(t,x)$
and is uniformly continuous on ${\mathbb T}_0^+\times B_1$. The equality $\dfrac{\partial f}{\partial x}(t,0)\equiv0$ is evident
since supports of $\Psi_{\ell}$ do not intersect $t$ axis. This implies $f\in{\cal F}$.
\hfill$\square$\medskip

\noindent\textbf{Remark 6.5.} Recall that in Example 4.6 we construct a linear system with constant matrix on a syndetic time scale which has
negative Lyapunov exponents but positive central upper exponent. Theorem 6.1 implies that this (asymptotically stable and even exponentially stable) 
system becomes unstable under arbitrarily small linear perturbation. Theorem 6.4 shows that a nonlinear system with exponentially stable first 
approximation can be unstable. Such examples can be found, for instance, in [26], but for non-regressive time scale systems.

\section{Stability and instability by first approximation}

First of all, we recall the time scale version of the Lyapunov theorem on asymptotic stability by first approximation, proved in 
[19], see also [27]. \medskip

\noindent\textbf{Definition 7.1. } Let $r>0$. We say that a continuous function $V(t,x):{\mathbb T}_0^+\times B_r\to {\mathbb R}$ is a \emph{strict Lyapunov function} for 
a time scale system
$$
x^\Delta=F(t,x), \qquad t\in {\mathbb T}_0^+, \quad x\in {\mathbb R}^n; \qquad F(t,0)\equiv 0, 
\eqno (7.1)$$
 if the following conditions are fulfilled for some $a>0$ and for all $t\in {\mathbb T}^+_a$, $x\in B_r$:
\begin{enumerate}
\item $V(t,x)\ge w_+(x)$, and $V(t,0)\equiv0$;
\item the trajectory $\Delta$-derivative of $V$ satisfies $V^{\Delta}(t,x)\le -w_-(x)$.
\end{enumerate}
Here $w_\pm(x):B_r \to {\mathbb R}$ are positive definite functions.\medskip 

\noindent\textbf{Remark 7.2.} Note that condition 2 means
\begin{eqnarray*}
\dfrac{\partial V}{\partial t}(t,x)+\dfrac{\partial V}{\partial x}(t,x)\cdot F(t,x)\le-w_-(x) & \quad \forall t\in {\mathbb T}^+_a\cap{\mathbb D};\\
V(\sigma(t),x+\mu(t)F(t,x))- V(t,x)\le-\mu(t)w_-(x) & \quad \forall t\in {\mathbb T}^+_a\cap{\mathbb S}.
 \end{eqnarray*}

\noindent\textbf{Theorem 7.3. (Lyapunov's Theorem) [19].} \emph{ If there is a strict Lyapunov function for system $(7.1)$, then the 
zero solution of this system is asymptotically stable.}\medskip

\noindent\textbf{Definition 7.4.} 
A constant matrix ${\cal A}$ is called \emph{strongly stable} with respect to $\mathbb T$ if its eigenvalues $\lambda_k$, $k=1,\dots, n$, satisfy inequality
 $\limsup\limits_{{\mathbb S}\ni t\to+\infty}\frac 1{\mu(t)}(|1+\mu(t)\lambda_k|^2-1)<0$.\medskip

 \noindent\textbf{Remark 7.5.} It is easy to see that if ${\cal A}$ is strongly stable then the following is true:
\begin{enumerate}
\item $\Re(\lambda_k)<0$, $k=1,\dots, n$;
\item time scale $\mathbb T$ is syndetic. 
\end{enumerate}
\medskip

\noindent\textbf{Theorem 7.6.} \emph{Suppose that the matrix ${\cal A}$ is strongly stable. Then there exists $\varepsilon>0$ such that for  any $r>0$ and any 
$f:{\mathbb T}^+_0\times B_r\to {\mathbb R}^n$ satisfying condition $(4.1)$, the solution $x=0$ of the system 
$$
x^\Delta= {\cal A}x+f(t,x)
\eqno (7.2)$$
is asymptotically stable.}\medskip

\noindent\textbf{Proof.} We show that, under the assumptions of theorem, there is a positive definite matrix ${\cal B}$ such that 
the quadratic form $V(x)=x^T {\cal B} x$ is a strict Lyapunov function for the system $(7.2)$. 

Making a non-degenerate transformation $x={\cal S}y$, we can reduce the first approximation system $(2.2)$ to the Jordan form
$$
y^\Delta={\cal J}y, 
\eqno (7.3)$$
where ${\cal J}={\cal S}^{-1}{\cal A}{\cal S}=\mbox{diag}({\cal J}_1,\ldots,{\cal J}_k)$ while for any $m=1,\ldots,k$
$${\cal J}_m=\begin{pmatrix} \lambda_m & 0 & 0 & \ldots & 0 & 0 \\
\delta & \lambda_m & 0 & \ldots & 0 & 0\\
 \ldots & \ldots & \ldots & \ldots & \ldots & \ldots\\
0 & 0 & 0 & \ldots & \delta & \lambda_m
\end{pmatrix}=:\lambda_m E +\delta I.
\eqno (7.4)$$
A parameter $\delta>0$ may be selected arbitrarily small.

System $(7.2)$ takes form 
$$
y^\Delta={\cal J}y+g(t,y),
$$
where $g(t,y)={\cal S}^{-1}f(t,{\cal S}y)$. The assumption $(4.1)$ implies 
$$|g(t,y)|\le |{\cal S}|\cdot|{\cal S}^{-1}|\cdot\varepsilon|y|\le C(\delta)\cdot\varepsilon|y|.
$$

First, we construct the desired quadratic form for the system $(7.3)$. It suffices to consider the system  
$$
z^\Delta={\cal J}_m z. 
\eqno (7.5)$$

We set $V(z)=|z|^2$ for $(7.5)$. Direct calculation of trajectory $\Delta$-derivative gives
$$V^{\Delta}=\begin{cases}
\displaystyle \frac {|1+\mu(t)\lambda_m|^2-1}{\mu(t)}|z|^2 + 2\delta\Re((1+\mu(t)\lambda_m)z\cdot I z)+\delta^2\mu(t)|I z|^2, & t\in {\mathbb S};\\
2\Re(\lambda_m) |z|^2 + 2\delta \Re(z\cdot I z), & t\in {\mathbb D}.
\end{cases}
$$
Taking into account Remark 7.5, we obtain $V^{\Delta}\le -\varkappa |z|^2$ with some $\varkappa>0$, if $\delta$ is sufficiently small and 
$t\in{\mathbb T}^+_a$ for $a$ sufficiently large.

For nonlinear system $(7.2)$ we set $V(y)=|y|^2$ and observe that $V^{\Delta}\le -\frac {\varkappa}2 |y|^2$ if $\varepsilon$ is sufficiently small.
\hfill$\square$\medskip

\noindent\textbf{Corollary 7.7. }\emph{If a matrix $\cal A$ is strongly stable with respect to $\mathbb T$ then the central upper exponent of the system
(2.2) is negative.}\medskip

\noindent\textbf{Proof.} By Theorem 7.6, there exists a transformation $x={\cal S}y$ such that the trajectory derivative of the Lyapunov function $V(y)=|y|^2$ 
satisfies $V^{\Delta}\le -\varkappa |y|^2$ with some $\varkappa>0$. This means that for any solution 
$\varphi(t)$ of linear system (2.2) we have
$$
|\varphi(t)|\le C e_{-\varkappa/2}(t,s) |\varphi(s)|,\qquad t>s, \quad t,s\in {\mathbb T},
$$
where $C$ is a positive constant depending on the transformation matrix $\cal S$. 

So, $u(t)\equiv-\varkappa/2$ is an upper function for system (2.2) and thus $\chi({\cal A})<0$.\hfill$\square$\medskip

Now we prove an analog of the famous Chetaev theorem on instability by first approximation (see [28] for the classical 
theorem and  [29] for the ``discrete'' one) for time scale systems.\medskip 

\noindent\textbf{Definition 7.8.} Let $r>0$. We say that a continuous function $V(t,x):{\mathbb T}_0^+\times B_r\to {\mathbb R}$ is a \emph{Chetaev function} 
for the system $(7.1)$ if the following conditions are fulfilled for some $a>0$ and for all $t\in {\mathbb T}^+_a$: 
\begin{enumerate}
\item $0\in {\mathfrak d} \Omega_t$, where $\Omega_t=\{x\in B_r: V(t,x)>0\}$; 
\item  $V$ is continuous at the origin uniformly with respect to $t$; 
\item the trajectory $\Delta$-derivative of $V$ satisfies $V^{\Delta}(t,x)\ge w(x)$, where $w(x):\Omega_t \to {\mathbb R}$ is a positive definite function 
(compare with condition 2 in Definition 7.1).
\end{enumerate}

\noindent\textbf{Theorem 7.9.}\emph{ If there is a Chetaev function for system $(7.1)$, then the zero solution of this system is unstable.}\medskip

\noindent\textbf{Proof.} 
Let $t_0\in {\mathbb T}^+_a$, and let $x_0\in \Omega_{t_0}$. Denote by $\widehat x(t)$ the solution of $(7.1)$ corresponding to initial conditions $x(t_0)=x_0$
 (Fig.\, 3). By condition 3, the function $V(t, \widehat x(t))$ increases while $\widehat x(t)\in B_r$. Moreover, the set 
$\{(t,x): x\in \Omega_t,  V(t,x)\ge V(t_0,x_0)\}$
is uniformly separated from zero by the condition 2, and therefore $(V(t,\widehat x(t)))^{\Delta}\ge b>0$. 
This means that $\widehat x(t)$ leaves the ball $B_r$, since otherwise $V(t,\widehat x(t))$ is unbounded. 
Since $x_0$ can be chosen arbitrarily close to zero by condition 1, the zero solution is unstable. \hfill$\square$\medskip

\begin{figure}[!ht]
\begin{center}
\includegraphics*[width=4in]{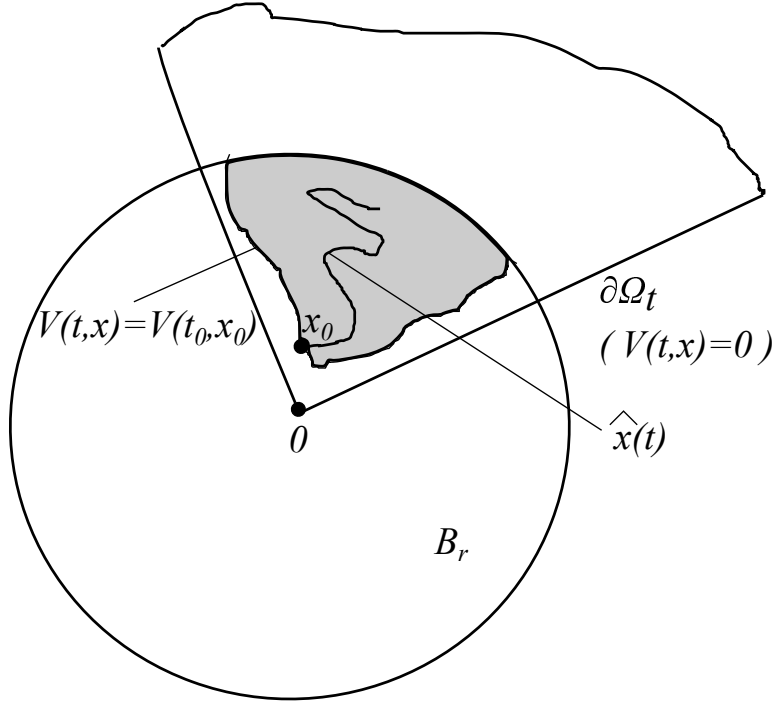}
\end{center}
\caption{Selection of $x_0$.}
\end{figure}

\noindent\textbf{Definition 7.10.} 
A constant matrix ${\cal A}$ is called \emph{strongly unstable} with respect to $\mathbb T$ if we can split its eigenvalues into two sets (the second one may be empty)
$$\lambda_k^{(1)}, \quad k=1,\dots, \ell;\qquad\qquad \lambda_j^{(2)}, \quad j=1,\dots, n-\ell;\qquad 1\le \ell\le n,
$$
such that the following inequalities are satisfied:
\begin{enumerate}
\item $\liminf\limits_{{\mathbb S}\ni t\to+\infty}\frac 1{\mu(t)}(|1+\mu(t)\lambda_k^{(1)}|^2-1)>0$, $k=1,\dots, \ell$;
\item $\liminf\limits_{{\mathbb S}\ni t\to+\infty}\frac 1{\mu(t)}(|1+\mu(t)\lambda_k^{(1)}|^2-|1+\mu(t)\lambda_j^{(2)}|^2)>0$, $k=1,\dots, \ell$, $j=1,\dots, n-\ell$;
\item if $\sup\mathbb D=+\infty$, we assume in addition that $\Re(\lambda_k^{(1)})>0$, $\Re(\lambda_k^{(1)})>\Re(\lambda_j^{(2)})$ for all 
$k=1,\dots, \ell$, $j=1,\dots, n-\ell$.
\end{enumerate}

\noindent\textbf{Remark 7.11.} Given a time scale $\mathbb T$, strongly stable and strongly unstable matrices form two non-intersecting classes. 
For time-invariant systems of ordinary differential equations (${\mathbb T}={\mathbb R}$) these matrices satisfy the assumptions $\max\limits_k\Re \lambda_k<0$
(Hurwitz matrices) and $\max\limits_k\Re \lambda_k>0$, respectively.\medskip

\noindent\textbf{Theorem 7.12.} \emph{ Let $\mathbb T$  be a syndetic time scale. Suppose that the matrix ${\cal A}$ is strongly unstable. 
Then there exists $\varepsilon>0$ such that for any $r>0$ and any $f:{\mathbb T}^+_0\times B_r\to {\mathbb R}^n$ satisfying condition $(4.1)$, 
the solution $x=0$ of the system $(7.2)$ is unstable.}\medskip

\noindent\textbf{Proof.} We show that, under the assumptions of theorem, there is a matrix ${\cal B}$ such that 
the quadratic form $V(x)=x^T {\cal B} x$ is a Chetaev function for the system $(7.2)$. 

As in the proof of Theorem 7.6, we can reduce the first approximation system $(2.2)$ to the Jordan form
$$
(y^{(1)})^\Delta={\cal J}^{(1)}y^{(1)};\qquad (y^{(2)})^\Delta={\cal J}^{(2)}y^{(2)}, 
\eqno (7.6)$$
where $y^{(1)}\in\mathbb R^{\ell}$, $y^{(2)}\in\mathbb R^{n-\ell}$, and, similarly to $(7.4)$,
$${\cal J}^{(1)}=\mbox{diag}(\lambda_1^{(1)} E +\delta I,\ldots,\lambda_m^{(1)} E +\delta I);\qquad 
{\cal J}^{(2)}=\mbox{diag}(\lambda_1^{(2)} E +\delta I,\ldots,\lambda_l^{(2)} E +\delta I).
$$
A parameter $\delta>0$ may be selected arbitrarily small.

We set $V(y)=|y^{(1)}|^2-|y^{(2)}|^2$ for $(7.6)$. Direct calculation of trajectory $\Delta$-derivative gives
$$V^{\Delta}\!=\begin{cases}
\displaystyle \sum_{i=1}^{\ell} \frac {|1+\mu(t)\lambda_i^{(1)}|^2-1}{\mu(t)}\ (y_i^{(1)})^2-\sum_{i=1}^{n-\ell} \frac {|1+\mu(t)\lambda_i^{(2)}|^2-1}{\mu(t)}\ (y_i^{(2)})^2
+ O(\delta)|y|^2, & t\in {\mathbb S};\\
\displaystyle 2\sum_{i=1}^{\ell} \Re(\lambda_i^{(1)}) (y_i^{(1)})^2 -2\sum_{i=1}^{n-\ell} \Re(\lambda_i^{(2)}) (y_i^{(2)})^2+O(\delta)|y|^2, & t\in {\mathbb D}.
\end{cases}
$$
By assumptions 1--3, we conclude that $V>0$ implies $V^{\Delta}\ge \varkappa |y|^2$ with some $\varkappa>0$, if $\delta$ is sufficiently small and 
$t\in{\mathbb T}^+_a$ for $a$ sufficiently large.

For nonlinear system $(7.2)$, similarly to Theorem 7.6, we obtain $V^{\Delta}\ge \frac {\varkappa}2 |y|^2$ if $\varepsilon$ is sufficiently small. \hfill$\square$

\section*{Appendix. Proof of Millionschikov's theorem (Theorem 6.2)}

We use the following relation (see [7, Page. 116, (8.8)]):
$$
\chi(\widetilde {\cal A})=\lim_{T\to\infty}\limsup_{k \to \infty}\dfrac{1}{kT} \sum_{i=0}^{k-1} \log | \Phi((i+1)T,iT)|,
\eqno (A.1)$$
where $\Phi$ is the Cauchy matrix for the system $(5.1)$.\medskip

We start with the main idea of the proof. Consider $T>0$ so that the value
$$\limsup_{k \to \infty}\dfrac{1}{kT} \sum_{i=0}^{k-1} \log | \Phi((i+1)T,iT)|$$
is close to $\chi(\widetilde {\cal A})$, see (A.1). Let $x_i$ $(i=0,1,2,\dots)$ be a unit vector such that
$$
| \Phi((i+1)T,iT)x_i| =| \Phi((i+1)T,iT)|,
\eqno (A.2)$$
and put $x_i(t)=\Phi(t,iT)x_i$.

It is $x_0(t)$ that has the fastest growth among solutions of $(5.1)$ on $[0,T]$. Without loss of generality, we may say that on $[T,2T]$, 
the solution $x_1(t)$ increases faster than $x_0(t)$.

We perturb system $(5.1)$  in the following way. First of all, we rotate the solution $x_0(t)$ in the plane 
$\langle x_0(t),x_1(t)\rangle$ by an angle $\delta>0$. Thus we obtain a function $y_0(t)$. This rotation can be done on a time segment of 
length $\ll T$. Then, for greater values of $t$, we set perturbation zero. Since $x_1(t)$ increases faster than $x_0(t),$ the angle between 
vectors $y_0(t)$ and $x_1(t)$ becomes less than $\delta$. This happens on a time period of length $\ll T$. Then we perturb system $(5.1)$ 
so that $y_0(t)$ becomes parallel to $x_1(t)$. Then we set the perturbation equal to zero up to $t=2T.$

Similarly, we consider segment $[2T,3T]$ and later ones. Finally, we obtain a solution $y_0(t)$ of the perturbed system that has Lyapunov 
exponent, close to $\chi(\widetilde {\cal A})$.\medskip

Now we proceed to the detailed proof.\medskip

\noindent\textbf{Step 1.} 
Given $\varepsilon>0$ and $\delta>0$, we fix a $T_0>1$ so that\footnote{It is sufficient to take $\varepsilon T_0/2$ instead of $\varepsilon T_0/4$ in formulae 
$(A.3)$--$(A.5)$. In fact, we need such selection of $T_0$ in Theorem 6.3 where we reproduce a part of the proof of Millionschikov's theorem.}
$$
\exp(\varepsilon T_0/4)\cdot \sin^2 \delta \geq 1.
\eqno (A.3) $$
Let triangles $\bigtriangleup ABC$ and $\bigtriangleup A_1B_1C_1$ be such that
$$
\dfrac{B_1C_1}{A_1C_1}:\dfrac{BC}{AC} \geq \exp(\varepsilon T_0/4);\quad
\sphericalangle A=\delta.
\eqno (A.4)$$
Then $(A.3)$, $(A.4)$ and Sine Theorem imply that
$$
\sin \sphericalangle B_1\leq \dfrac{\sin \sphericalangle B_1}{\sin \sphericalangle A_1}= \dfrac{A_1C_1}{B_1C_1} 
\leq\exp(-\varepsilon T_0/4)\cdot \dfrac{1}{\sin \delta }\leq \sin \delta.
\eqno (A.5)$$
Since $A_1C_1/B_1C_1 \leq 1$ we have $\sphericalangle B_1\leq \sphericalangle A_1,$ and, consequently,
$\sphericalangle B_1\leq \pi/2.$
Therefore, $(A.5)$ implies $\sphericalangle B_1\leq \delta$. \medskip

\noindent\textbf{Step 2.} 
Fix $T>0$ such that $m=T/T_0\in {\mathbb N}$, $(2a+(2a+1)\delta)/m<\varepsilon/8$ and
$$
\limsup_{k \to \infty} \dfrac{1}{kT} \sum_{i=0}^{k-1} \log | \Phi((i+1)T,iT)| > \chi(\widetilde {\cal A}) -\dfrac{\varepsilon }{4}.
\eqno (A.6) $$

\noindent\textbf{Step 3.} Take a unit vector $x_i\: (i=0,1,2,\dots)$ such that $(A.2)$ is satisfied.
Let
$$
x_i(t)=\Phi(t,iT)x_i
\eqno (A.7) $$
be solutions of $(5.1)$.

Set $\widehat {\cal B}(t)=0$ for $0 \leq t \leq T$. Suppose that
$$\dfrac{| x_1(2T)|}{| x_1(T)|}:\dfrac{| x_0(2T)|}{| x_0(T)|} > \exp(3\varepsilon T/4).
$$
(if this is wrong, we set $\widehat {\cal B}(t)=0$ for $T<t \leq 2T$).

Divide the segment $[T,2T]$ to $m$ segments of length $T_0$:
$$Q_l=[\tau_l,\tau_{l+1}]=[T+(l-1)T_0,T+lT_0]\quad (l=1,2,\dots,m).
$$ 
Let $Q_{j}$ be the first of segments  $Q_1,Q_2,\dots ,Q_{m-1}$ 
where
$$
\dfrac{| x_1(\tau_{l+1})|}{| x_1(\tau_l)|}:\dfrac{| x_0(\tau_{l+1})|}{| x_0(\tau_l)|} \geq
\exp(\varepsilon T_0/2).
\eqno (A.8)$$
So, $\tau_j <\tau_{j+1} <\tau_{j+2}$ are ends of segments  $Q_{j}$ and $Q_{j+1}$ (Fig.\, 4).

\begin{figure}[!ht]
\begin{center}
\includegraphics*[width=4in]{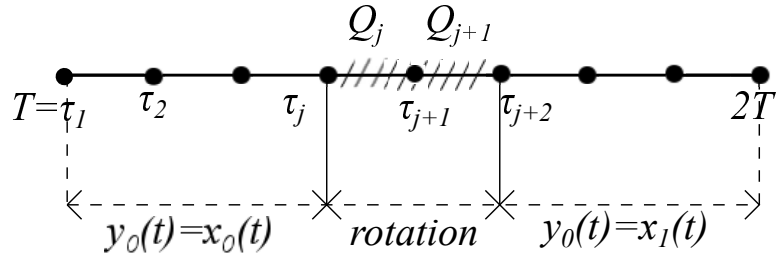}
\end{center}
\caption{Segments $Q_l$.}
\end{figure}

Note that the number of values of $l$ for which $(A.8)$ is satisfied is not less than $2$. Indeed, otherwise
$$\dfrac{| x_1(2T)|}{| x_1(T)|}:\dfrac{| x_0(2T)|}{| x_0(T)|} \le \exp((m-1)\varepsilon T_0/2)\cdot 
\dfrac{| x_1(\tau_{l+1})|}{| x_1(\tau_l)|}:\dfrac{| x_0(\tau_{l+1})|}{| x_0(\tau_l)|}.
$$
Since for any nonzero solution $x(t)$ of $(5.1)$ we have 
$$
\exp(-aT_0)\le \dfrac{|x(\tau_{l+1})|}{|x(\tau_{l})|}\le \exp(aT_0),
$$
(we recall that $a=\|\widetilde {\cal A}\|_{{\mathbb L}^\infty}$), this implies
$$\dfrac{| x_1(2T)|}{| x_1(T)|}:\dfrac{| x_0(2T)|}{| x_0(T)|} \le \exp(\varepsilon T/2)\cdot \exp(2aT_0)\le \exp(5\varepsilon T/8),
$$
a contradiction. Therefore, $\tau_{j+1}<2T$.

Define the perturbation $\widehat {\cal B}(t)$ for $T<t \leq 2T$ in the following way.

\noindent\textbf{A.} If $t \notin [\tau_{j},\tau_{j}+1]\bigcup [\tau_{j+1},\tau_{j+1}+1]$ we set 
$\widehat {\cal B}(t)=0$.

\noindent\textbf{B.} For $t \in [\tau_{j},\tau_{j}+1]$ we set
$$
\widehat {\cal B}(t)=U_\delta^{-1} (t)\widetilde {\cal A}(t)U_\delta(t)-U_\delta^{-1}(t)\dot U_\delta(t)-\widetilde {\cal A}(t),
\eqno (A.9)$$
where $U_\delta(t)$ is an orthogonal matrix such that
$$
U_\delta(\tau_{j})=E_n, \qquad
| \dot U_\delta (t) | \leq \delta.
\eqno (A.10)$$

Namely, we define $U_\delta(t)$ as a rotation in the plane $\langle x_0(\tau_{j}+1), x_1(\tau_{j}+1)\rangle$
in the direction from $x_0(\tau_{j}+1)$ to $x_1(\tau_{j}+1)$ with the speed not greater then $\delta$. 
From $(A.9)$, $(A.10)$ and orthogonality of $U_\delta$ we deduce inequality $(6.3)$.

By construction, there exist $\alpha_1 \geq 0, \alpha_2>0$ such that 
$$
U_\delta^{-1}(\tau_{j}+1)x_0(\tau_{j}+1)=
\alpha_1 x_0(\tau_j+1)+\alpha_2 x_1(\tau_{j}+1),
\eqno (A.11)$$
and
$$
\sphericalangle {(x_0(\tau_{j}+1),y_0(\tau_{j}+1))}=\delta.
\eqno (A.12)$$

\noindent\textbf{C.} Due to relations $(A.11)$, $(A.12)$, $(A.8)$ and to statements of {\bf Step 1}, we have
$$\sphericalangle (\alpha_1 x_0(\tau_{j+1}+1)+\alpha_2 x_1(\tau_{j+1}+1), x_1(\tau_{j+1}+1)) \leq \delta
$$
(here $\alpha_1$ and $\alpha_2$ are defined by $(A.11)$).

For $t \in [\tau_{j+1},\tau_{j+1}+1]$  we take $\widehat {\cal B}(t)$ that satisfies $(A.9)$ and $(A.10)$ (with $\tau_{j}$ replaced by $\tau_{j+1}$). 
Instead of inequalities $(A.11)$ and $(A.12)$ we demand that
$$U_\delta^{-1}(\tau_{j+1}+1)(\alpha_1 x_0(\tau_{j+1}+1)+\alpha_2 x_1(\tau_{j+1}+1))=\beta x_1(\tau_{j+1}+1)
$$ 
for some $\beta>0$ (Fig.\,5).

\begin{figure}[!ht]
\begin{center}
\includegraphics*[width=4in]{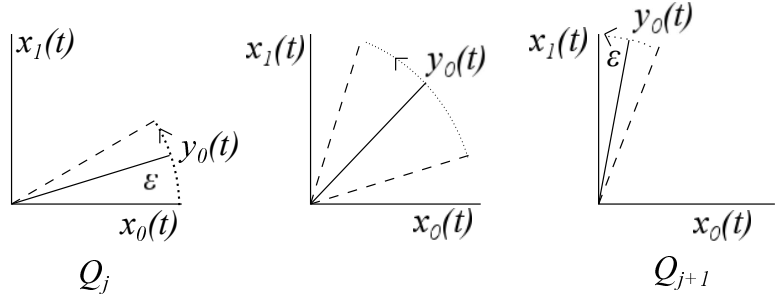}
\end{center}
\caption{Millionschikov's rotations.}
\end{figure}

Observe that since $x_0(t)$ is a solution of $(5.1)$, the function
$$y_0(t)=\left\{
\begin{array}{rl}
x_0(t) &\mbox{\quad if \quad} T\le t\le\tau_{j}\\
U_\delta^{-1}(t)x_0(t) &\mbox{\quad if \quad} \tau_{j}\le t\le\tau_{j}+1,\\
\alpha_1 x_0(t)+\alpha_2 x_1(t) &\mbox{\quad if \quad} \tau_{j}+1\le t\le\tau_{j+1},\\
U_\delta^{-1}(t)(\alpha_1 x_0(t)+\alpha_2 x_1(t)) &\mbox{\quad if \quad} \tau_{j+1}\le t\le\tau_{j+1}+1,\\
\beta x_1(\tau_{j+1}+1) &\mbox{\quad if \quad} \tau_{j+1}+1\le t\le 2T,\\
\end{array}\right.$$
is a solution of system $(5.3)$ with constructed matrix $\widehat {\cal B}(t)$.\medskip

\noindent\textbf{Step 4.} We construct the perturbation $\widehat {\cal B}(t)$ on segments $[iT,(i+1)T]$, $i=2,3,\dots$, basing on solution $x_i(t)$
similarly to what we have done above.\medskip

\noindent\textbf{Step 5.} Consider the constructed solution $y_0(t)$ of the system $(5.3)$. We claim that $y_0(t)$ has the Lyapunov exponent greater than 
$\chi(\widetilde {\cal A})-\varepsilon$. Indeed, due to $(A.2)$, $(A.6)$ and $(A.7)$ it suffices to prove that for any $i=0,1,2,\dots$
$$
\dfrac{| y_0((i+1)T)|}{| y_0(iT)|} \geq
\dfrac{| x_i((i+1)T)|}{| x_i(iT)|}\exp\left(-\dfrac{3\varepsilon T}{4}\right).
$$
It follows from construction of $y_0(t)$ that for any fixed $i$ the number of $l$ such that inequality
$$
\dfrac{| y_0(iT+(l+1)T_0)|}{| y_0(iT+lT_0)|} \geq
\dfrac{| x_i(iT+(l+1)T_0)|}{| x_i(iT+lT_0)|}\exp\left(-\dfrac{\varepsilon T_0}{2}\right)
\eqno (A.13) $$
is not fulfilled, does not exceed 2 (for $i=1$ this might be only segments $Q_j$ and $Q_{j+1}$, see Fig.~4). If $(A.13)$ is not satisfied, 
we use the inequality
$$
\dfrac{| y_0(iT+(l+1)T_0)|}{| y_0(iT+lT_0)|} \geq
\dfrac{| x_i(iT+(l+1)T_0)|}{| x_i(iT+lT_0)|}\exp(-(2a+(2a+1)\delta)T_0).
\eqno (A.14)$$

Multiplying inequalities $(A.13)$ and $(A.14)$ corresponding to $l=0,1,\dots,m-1$,
we obtain
\begin{eqnarray*}
\dfrac{| y_0((i+1)T)|}{| y_0(iT)|} &>& \dfrac{| x_i((i+1)T)|}{| x_i(iT)|} 
\exp\left(-\left(\dfrac{\varepsilon}{2} + \dfrac{2(2a+(2a+1)\delta)T_0}{T}\right)T\right) \\
&\ge& \dfrac{| x_i((i+1)T)|}{| x_i(iT)|}\exp\left(-\dfrac{3\varepsilon T}{4} \right)
\end{eqnarray*}
(the last inequality holds by the choice of $T$ in {\bf Step 2}). This completes the proof. \hfill$\square$\medskip

\noindent{\bf Corollary A.1.} {\em The perturbation $\widehat {\cal B}(t)$ may be taken continuous.}\medskip

\noindent\textbf{Proof.} It follows from the proof that $\widehat {\cal B}(t)$ is piecewise continuous i.e. has finitely many discontinuity points 
on bounded subsets of $\mathbb R$. So, we may construct a  continuous matrix $\widehat {\cal B}_1(t)$ such that 
$|\widehat {\cal B}_1(t)|\le(2a+1)\delta$, and
$$M=\{t: \widehat {\cal B}(t)\neq \widehat {\cal B}_1(t)\}=\bigcup\limits_{j=1}^{\infty}\Delta_j,
$$
where the length of intervals $\Delta_j$ can be chosen arbitrarily small.

Consider the system
$$
\dot x=(\widetilde {\cal A}(t)+\widehat {\cal B}_1(t))x.
\eqno  (A.15) $$
Let $\Psi(t)$ and $\Xi(t)$ be fundamental matrices of $(A.15)$ and $(5.3)$ respectively, so that $\Psi(0)=\Xi(0)=E_n$. Then
$$\Xi(t)=\Psi(t)+\Psi(t)\int_0^t\Psi^{-1}(\tau)(\widehat {\cal B}(\tau)-\widehat {\cal B}_1(\tau))\Xi(\tau)\, d\tau,$$
and thus
$$
|\Xi(t)| \leq |\Psi(t)| +|\Psi(t)| \int_0^t | \Psi^{-1}(\tau)||\widehat {\cal B}(\tau)-\widehat {\cal B}_1(\tau)||\Xi(\tau)|\, d\tau.
\eqno (A.16) $$

Denote $u(t)=|\Xi(t)| / | \Psi(t)|,$ 
$v(t)=| \Psi^{-1}(t)|| \Psi(t)||\widehat {\cal B}(t)-\widehat {\cal B}_1(t)|$.
Dividing both parts of $(A.16)$ by $| \Psi(t)|,$ we obtain
$$ u(t)\leq 1+ \int_0^t u(\tau) v(\tau)\, d \tau, 
$$ 
which implies by the Gr\"{o}nwall--Bellmann lemma
$$\dfrac{|\Xi(t)|}{|\Psi(t)|}=u(t) \leq \exp\left(\int_0^t v(\tau)\, d \tau\right)\le
\exp\Big((2a+1)\delta\sum\limits_{j=1}^{\infty}|\Delta_j|\cdot\sup\limits_{s\in\Delta_j}\exp(2(a+(2a+1)\delta)s)\Big).
$$
The last expression can be made arbitrarily close to $1$, and the statement follows. \hfill$\square$

\bigskip

\noindent\textbf{Acknowledgements.} The first author was partially supported by RFBR grants 
14-01-00202 and 15-01-03797-a, by St.Petersburg State University under Thematic Plan 6.38.223.2014, and by the Fulbright Program. 
The second author was supported by RFBR grant 15-01-07650 and by St.Petersburg State University grant 6.38.670.2013. 

Authors are grateful to Prof. Martin Bohner, Prof. Anatoly Martynyuk and Prof. Andrejs Reinfelds for their attention to our research and 
for their precious advices and remarks. We also thank Prof. Vitaly Slyn'ko who provided us with reference [26], and anonymous referee for useful comments.


\begin{thebibliography}{99}

\bibitem{krynaz1} B. Aulbach,\, S. Hilger, \emph{Linear Dynamic Processes with Inhomogenous Time Scale}, In Nonlinear Dynamics and Quantum Dynamical Systems (Gaussig, 1990), 
volume 59 of Math. Res., pages 9--20. Akademie Verlag, Berlin, 1990.

\bibitem{krynaz2} A.\, Reinfelds, L.\, Sermone, \emph{Stability of Impulsive Differential Systems}, Abstr. Appl. Anal.
\textbf{2013} (2013), Article ID 253647, 11 pages.

\bibitem{krynaz3} M.\, Bohner, \emph{Some Oscillation Criteria for First Order Delay Dynamic Equations}, Far East J. Appl. Math. \textbf{18:3} (2005), 289--304.

\bibitem{krynaz4} P.\,E.\, Kloeden, A.\, Zmorzynska, \emph{Lyapunov Functions for Linear Nonautonomous Dynamical Equations on Time Scales},  Adv. Differ. Equ, Article ID69106, 
\textbf{2006} (2006), 1--10.

\bibitem{krynaz5} M.\, Bohner, D.\,A.\, Lutz, \emph{Asymptotic Behavior of Dynamic Equations on Time Scales}, J. Differ. Equations Appl.,  \textbf{7:1} (2001), 21--50.

\bibitem{krynaz6} S.\, Bodine, D.\,A.\, Lutz, \emph{Exponential Functions on Time Scales: Their Asymptotic
Behavior and Calculation}, Dynam. Systems Appl., \textbf{12} (2003), 23--43.

\bibitem{krynaz7} B.\,F.\, Bylov, R.\,E.\, Vinograd, D.\,M.\, Grobman, V.\,V.\, Nemytskii, \emph{Teoriya pokazatelei Lyapunova i ee prilozheniya k voprosam ustoichivosti} 
(Theory of Lyapunov Exponents and its Application to Problems of Stability), Moscow: Nauka, 1966, 576 p. (in Russian).

\bibitem{krynaz8} M.\, Bohner, A.\,A.\, Martynyuk, \emph{Elements of Stability Theory of A.M. Liapunov for Dynamic Equations on Time Scales}, Nonlinear 
Dynamics and Systems Theory, \textbf{7:3} (2007), 225--251.

\bibitem{krynaz9} N.\, H. Du, L.\, H. Tien, \emph{On the Exponential Stability of Dynamic Equations on Time Scales}, J. Math. Anal. Appl. \textbf{331} (2007), 1159--1174.

\bibitem{krynaz10} J.\, Hoffacker,  C.\,C.\, Tisdell, \emph {Stability and Instability for Dynamic Equations on Time Scales}, Comput. Math. Appl., \textbf{49:9--10} 
(2005), 1327--1334.

\bibitem{krynaz11} A.\, A.\, Martynyuk, \emph{On the Exponential Stability of a Dynamical System on a Time Scale}, Dokl. Akad. Nauk. \textbf{421} (2008), 312--317.

\bibitem{krynaz12} T.\, Gard, J.\, Hoffacker, \emph{Asymptotic Behavior of Natural Growth on Time Scales}, Dynam. Systems Appl., \textbf{12:1--2} (2003), 131--148. 

\bibitem{krynaz13} G.\, Hovhannisyan, \emph{Asymptotic Stability for Dynamic Equations on Time Scales}, Adv. Difference Equ., \textbf{2006} (2006), Article ID 18157, 1--17.

\bibitem{krynaz14} G.\, Hovhannisyan, \emph{Asymptotic Stability for 2x2 Linear Dynamic Systems on Time Scales},
International Journal of Difference Equations, \textbf{2:1}  (2007), 105--121.

\bibitem{krynaz15} W.\,N.\, Li, \emph{Some Pachpatte Type Inequalities on Time Scales}, Computers and Mathematics with Applications, \textbf{57} (2009), 275--282.

\bibitem{krynaz16} D.\,B.\, Pachpatte, \emph{Explicit Estimates on Integral Inequalities with Time Scale}, J. Inequal. Pure Appl. Math. \textbf{7:4}, Article 143, (2006), 1--8.

\bibitem{krynaz17} A.\, I.\, Bobenko, Yu.\, B.\, Suris, \emph{Discrete Differential Geometry. Integrable Structure.}
Graduate Studies in Mathematics , Vol. 98. AMS, 2008. xxiv+404 p.

\bibitem{krynaz18} S.\,G.\, Kryzhevich, \emph{The Relation Between Central Exponents of Linear Systems of Ordinary Differential Equations and Stability of the Excited Systems},
Differential Equations, \textbf{36:10} (2000), 1430--1431.

\bibitem{krynaz19} M.\, Bohner, A.\, Peterson, \emph{Dynamic Equations on Time Scales. An Introduction with Applications,} Birkh\"{a}user Boston Inc., Boston, MA, 2001. 

\bibitem{krynaz20} G.\,Sh.\, Guseinov, \emph{Integration on Time Scales}, J. Math. Anal. Appl. \textbf{285} (2003), 107--127.

\bibitem{krynaz21} M.\, Bohner, A.\, Peterson, \emph{Advances in Dynamic Equations on Time Scales.} Birkh\"{a}user Boston Inc., Boston, MA, 2003.

\bibitem{krynaz22} S.\,K.\, Choi, D.\,M.\, Im, N.\, Koo, \emph{Stability of Linear Dynamic Systems on Time Scales}, Advances in Difference Equations,  Article ID 670203 
\textbf{2008} (2008), 1--12.

\bibitem{krynaz23} J.\,J.\, DaCunha, \emph{Stability for Time Varying Linear Dynamic Systems on Time Scales}, J. Comput. Appl. Math., \textbf{176:2} (2005): 381--410.

\bibitem{krynaz24} O.\, Perron, \emph{\"Uber Stabilit\"at und Asymptotisches Verhalten der Integrale von
Differentialgleichungssystemen}, Math. Zeitschrift. \textbf{29} (1928), 129--160 (In German).

\bibitem{krynaz25} V.\,M.\,  Millionschikov, \emph{A Proof of Attainability for Central Exponents of Linear Systems}, Siberian Mathematical Journal \textbf{10:1} (1969), 
99--104.

\bibitem{krynaz26} C.\, P\"otzsche, S.\, Siegmund, F.\, Wirth, \emph{A Spectral Characterization of Exponential Stability for Linear Time-Invariant Systems on Time Scales}, 
Discrete Contin. Dyn. Syst. \textbf{9} (2003), 1223--1241.

\bibitem{krynaz27} A.\, BenAbdallah, M.\, Dlala, M.\,A.\, Hammami, \emph{A New Lyapunov Function for Stability of 
Time-varying Nonlinear Perturbed Systems}, Syst. Control Lett. \textbf{56} (2007), 179--187.

\bibitem{krynaz28} N.\,G.\, Chetaev, \emph{The Stability of Motion.} New York: Pergamon Press. 1961.

\bibitem{krynaz29} C.\, Carcamo, C.\, Vidal, \emph{The Chetaev Theorem for Ordinary Difference Equations}, Proyecciones Journal of Mathematics \textbf{31:4} (2012), 391--402.

\end{thebibliography}
\end{document}